\documentclass[reqno, a4paper,12pt]{amsart}

\usepackage[capposition=top]{floatrow}
\usepackage{euscript,eufrak,verbatim}
\usepackage{graphicx}
\usepackage[usenames]{color}
\usepackage[colorlinks,linkcolor=red,anchorcolor=blue,citecolor=blue]{hyperref}
\usepackage{amsmath, mathtools}
\usepackage{mathtools}
\usepackage{amsthm}
\usepackage[all]{xy}
\usepackage{amssymb}
\usepackage{bm}

\usepackage[all]{xy}

\usepackage{mathrsfs}
\usepackage{amscd}

\usepackage{enumitem}

\makeatletter
\numberwithin{equation}{section}

\theoremstyle{plain}
\newtheorem{main theorem}{Main Theorem}
\newtheorem{theorem}{Theorem}[section]
\newtheorem{lemma}[theorem]{Lemma}

\newtheorem{corollary}[theorem]{Corollary}
\newtheorem{proposition}[theorem]{Proposition}
\newtheorem{claim}[theorem]{Claim}

\theoremstyle{definition}
\newtheorem{definition}[theorem]{Definition}
\newtheorem{remark}[theorem]{Remark}

\newcommand{\spa}{\hspace{1pt}}
\newcommand{\vep}{\varepsilon}
\newcommand{\flo}{\mathscr}
\newcommand{\diam}{\mathrm{diam}}
\newcommand{\tpi}{\widetilde{\pi}}

\begin{document}
\title[VARIATIONAL PRINCIPLE OF HIGHER DIMENSION WEIGHTED PRESSURE]{Variational principle of higher dimension weighted pressure for amenable group actions}

\author[Zhengyu Yin, Zubiao Xiao]{Zhengyu Yin, Zubiao Xiao}



\maketitle

\begin{abstract}
Let $r\geq 2$ and $(X_i,G)$ $(i=1,\cdots,r)$ be topological dynamical systems with $G$ being an infinite discrete amenable group. Suppose that $\pi_i:(X_i,G)\to (X_{i+1},G)$ are factor maps and $0\leq w_i\leq 1$. In this article, for $f\in C(X_1)$, we introduce the weighted topological pressure $P^{\textbf{a}}(f,G)$ for higher dimensions (not only for $r=2$) of amenable group actions. By using measure-theoretical theory, we establish a variational principle as
\begin{align*}
    P^{\textbf{a}}(f,G)=\sup_{\mu\in \mathcal{M}^G(X_1)}\Big(\sum_{i=1}^rw_ih_{\mu_i}(X_i,G)+w_1\int_{X_1}fd\mu\Big),
\end{align*}
where $\mu_i=\pi_{i-1}\circ\cdots\circ\pi_{1}\mu$ is the induced $G$-invariant measure on $X_{i}$.
\end{abstract}

\section{Introduction} \label{section: introduction}
Let $X$ be a compact metric space, $T$ a continuous self-map of $X$, and $\mu$ a $T$-invariant Borel probability measure on $X$. To distinguish the topological dynamical system $(X,T)$ and measure-theoretical system $(X,T,\mu)$, Kolmogorov \cite{ko} introduced the measure theoretical entropy $h_\mu(X, T)$ in 1958. Later in 1965, by measuring the complexity of open covers of $X$, Adler, Konheim, and McAndrew \cite{akm} defined the topological entropy for $(X,T)$. Although the two entropies are quite different in definitions, there is a variational principle to connect them. Goodman \cite{g} in 1969 and Goodwyn \cite{gw} in 1970 published the first variational principle in two separate papers.

As a generalization of topological entropy, topological pressure was first introduced by Ruelle \cite{ru}, and further studied by Walters \cite{pw1}, who developed a variational principle, that is,
\begin{align*}
    P(T,f)=\sup_{\mu\in \mathcal{M}(X,T)}\left\{h_{\mu}(X,T)+\int_X fd\mu\right\}\text{ for all $f\in C(X)$}.
\end{align*}
where $\mathcal{M}(X,T)$ is the set of $T$-invariant Boral probability measures on $X$.

Let $(X,T)$ and $(Y,S)$ be two dynamical systems. A map $\pi: X\to Y$ is said to be a \textit{factor} map if $\pi$ is a continuous and surjective map such that $S\circ\pi(x)=\pi\circ T(x)$ for all $x\in X$. Let $\boldsymbol{a}=(a_1,\cdots,a_k)\in \mathbb{R}^k$ with $a_1>0$ and $a_i\geq 0$ for $2\leq i\leq k$. Let $(X_i,T_i)$ be dynamical systems and $\pi:X_{i}\to X_{i+1}$ factor maps. Using Bowen's and Pesin-Pitskel's approaches, Feng and Huang in \cite{fh} introduced the $\boldsymbol{a}$-weighted topological pressure $\mathcal{P}^{\boldsymbol{a}}(T_1,f)$ for every $f\in C(X_1)$ and proved a variational principle for it, that is,
\begin{align*}
    \mathcal{P}^{\boldsymbol{a}}(T_1,f)=\sup_{\mu\in \mathcal{M}(X_1,T_1)}\left\{\sum_{i=1}^ra_ih_{\mu_i}(X_i,T_i)+\int_{X_1}fd\mu\right\},
\end{align*}
where $\mu_i=\pi_{i-1}\circ\cdots\circ\pi_{1}\mu$ is the induced $T_i$-invariant measure on $X_{i}$. In 2022, M. Tsukamoto \cite{mt} gave a new definition $P^{(a_1,a_2)}(X,f)$ to weighted topological pressure for the case of two dynamical systems $(X, T)$ and $(Y, S)$. He proved that for a factor map $\pi: X\to Y$ and any $a_1,a_2\geq 0$ with $a_1+a_2=1$,
\begin{align*}
    P^{(a_1,a_2)}(X,f)=\sup_{\mu\in \mathcal{M}(X,T)}\left\{a_1h_\mu(X,T)+a_2h_{\pi^*\mu}(Y,S)+a_1\int_{X_1}fd\mu\right\}, \text{ for all $f\in C(X)$},
\end{align*}
where $\pi^*\mu=\mu\circ \pi^{-1}$. Moreover, when $f=0$, he explained the importance of the equivalence of the two weighted topological entropy (Feng, Huang and M. Tsukamoto's different definitions, see \cite[Example 1.6]{mt}). Recently, following the approach of M. Tsukamoto, the authors J. Yang, E. Chen etc \cite{YCYY}. extended the $(a_1,a_2)$-weighted topological pressure of factor map $\pi: X\to Y$ to the action of amenable groups.

Let $(X_i,T_i)$ $(i=1,\cdots,r \text{ and }r\geq 2)$ be topological dynamical systems and $\pi_i:X_i\to X_{i+1}$ factor maps. N. Alibabaei \cite{al} redefined the weighted topological weighted pressure and entropy of M. Tsukamoto for higher dimensions and proved a variational principle. He also calculated the Hausdorff dimension and entropy for some special sets by using the coincidence of two definitions of Feng, Huang and M. Tsukamoto.

In this paper, we aim to generalize M. Tsukamoto's definition to higher dimensional (not less than $2$) factor maps for the action of amenable groups and obtain a variational principle. In section 2, we recall some fundamental notations and define the weighted topological pressure and entropy for factor maps and discuss the method of zero-dimensional principal extension for the action of amenable groups. In section 3 and section 4, we prove the variational principle between weighted topological pressure and weighted measure-theoretical pressure. Finally, as a simple application, we show that weighted pressure determines weighted measure-theoretic entropy for an amenable group action.
\section{Weighted topological pressure and entropy} \label{section: weighted topological pressure}

\subsection{Amenable group and convergence lemma}
Let $G$ be a countable discrete infinite group and denote $\mathcal{F}(G)$ the class of all finite non-empty subsets of $G$. Recall that a group $G$ is \textit{amenable} if it admits a sequence of non-empty finite subsets $\{F_n:F_n\in \mathcal{F}(G), n\in \mathbb{N}\}$ such that for any  $K\in \mathcal{F}(G)$, we have
\begin{center}
    $\underset{n\to \infty}{\lim}\,\displaystyle\frac{|B(F_n,K)|}{|F_n|}=0$,
\end{center}
where $B(F_n,K)=\{g\in G:Kg\cap F_n\ne \emptyset$ and $Kg\cap (G\backslash F_n)\ne\emptyset\}$ and $|\cdot|$ means the size of a set. We call the sequence $\{F_n:F_n\in \mathcal{F}(G), n\in\mathbb{N}\}$ an \textit{F$\o$lner sequence} of $G$. Throughout the paper, we assume $G$ is an infinite countable amenable group.

Let $\phi:\mathcal{F}(G)\to \mathbb{R}$ be a function, we say that $\phi$ is
\begin{enumerate}
    \item [(1)]\textbf{non-negative}, if $\phi(F)\geq 0$, for all $F\in \mathcal{F}(G)$.
    \item [(2)]  \textbf{monotone}, $\phi(E)\leq\phi(F)$, for all $E,F\in \mathcal{F}(G)$ with $E\subset F$.
    \item[(3)] \textbf{$G$-invariant} $\phi(Fg)=\phi(F)$, for all $F\in \mathcal{F}(G)$.
    \item[(4)] \textbf{sub-additive}, $\phi(E\cup F)\leq\phi(E)+\phi(F)$, for all $E,F\in \mathcal{F}(G)$ with $E\cap F=\emptyset$.
\end{enumerate}

The following Theorem for invariant sub-additive functions on finite subsets of amenable groups is due to Ornstein and Weiss (see \cite{lw},\cite{ow} for more details). It plays a central role in the definition of entropy, mean topological dimension theories, and so on.
\begin{theorem}\label{thm2.1}
Let $G$ be a countable amenable group. Let $\phi:\mathcal{F}(G)\to \mathbb{R}$ be a map satisfying the properties (1), (2), (3), and (4) above.
Then there is a real number $\lambda\geq 0$ such that \[\underset{n\to\infty}{\lim}\,\displaystyle\frac{\phi(F_n)}{|F_n|}=\lambda\] exists for all F$\o$lner sequence $\{F_n\}$ of $G$. Moreover, the limit is independent of the choice of F$\o$lner sequence in $G$.
\end{theorem}

\subsection{Weighted pressure and entropy}
We call $(X,G)$ a \textit{$G$-space} if $G$ acts on the topological space $X$ as a topological group with $ex=x$ and $g_1(g_2x)=(g_1g_2)x$ for the identity $e\in G$, any $g_1,g_2\in G$ and $x\in X$. Denote $\mathcal{M}(X)$ by the set of Borel probability measures and $\mathcal{M}^G(X)$ by the set of $G$-invariant Boral probability measures on $X$. Both $\mathcal{M}(X)$ and $\mathcal{M}^G(X)$ are endowed with weak*-topology. For any $\mu\in \mathcal{M}^G(X)$, we denote $h_\mu(X,G)$ the \textit{measure-theoretical entropy} of the measurable $G$-system $(X,G,\mu)$ (See the definition in section 4).

Let $(X,G)$ and $(Y,G)$ be two $G$-systems. A map $\pi: X\to Y$ is said to be a \textit{factor} map if $\pi$ is a continuous and surjective such that $g\circ\pi(x)=\pi\circ g(x)$ for all $x\in X$ and $g\in G$. We also say that $(X, G)$ is an extension of $(Y, G)$ or $(Y, G)$ is a factor of $(X, G)$.

Let $\pi:X\to Y$ be a factor map. For any $\mu\in \mathcal{M}^G(X)$, we put $\pi\mu\in \mathcal{M}^G(Y)$ as the induced $G$-invariant Borel measure on $Y$, that is, $\pi\mu(A)=\mu(\pi^{-1}(A))$ for all Borel set $A\subset Y$.

In this part, we introduce the definition of higher dimensional weighted topological pressure of factor maps for amenable group actions.

For a $G$-system $(X, G)$, we call a collection  $\mathcal{U}$ of subsets \textit{a cover} of $X$, if $\underset{U\in\mathcal{U}}{\bigcup}U=X$. We denote by $\mathcal{C}^o_X$ (respectively $\mathcal{C}_X$, $\mathcal{P}_X$) the class of finite open covers (respectively finite covers, finite partition) of $X$. Given $\mathcal{U},\mathcal{V}\in \mathcal{C}_X$, we write $\mathcal{V}\preceq\mathcal{U}$ if for any $U\in \mathcal{U}$, there is some $V\in \mathcal{V}$ such that $U\subset V$, we also set $\mathcal{V}\vee \mathcal{U}=\{V\cap U:V\in\mathcal{V}, U\in \mathcal{U}\}$.
Let $F\in \mathcal{F}(G)$, $\mathcal{U}\in \mathcal{C}_X$, put $\mathcal{U}_F=\bigvee_{g\in F}g^{-1}\mathcal{U}=\{\bigcap_{g\in F}g^{-1}U_g:U_g\in \mathcal{U}\}$.

Let $d$ be a metric on $X$ and $F\in\mathcal{F}(G)$. We define a metric $d_F$ on $X$ by
\begin{equation}\label{el2.1}
    d_F(x_1,x_2)=\underset{g\in F}{\max}\,d(gx_1,gx_2).
\end{equation}

Consider dynamical systems $(X_i,G)$ ($i=1,2,\cdots,r$) and factor maps $\pi_i:X_i\to X_{i+1}$ $(i=1,2,\cdots,r)$:
\begin{equation*}
    (X_1,G)\overset{\pi_1}{\longrightarrow}(X_2,G)\overset{\pi_2}{\longrightarrow}\cdots\overset{\pi_{r-1}}{\longrightarrow}(X_r,G).
\end{equation*}

Write
\begin{align*}
    \pi^{(0)}=id_{X_1}:X_1\to X_1,
\end{align*}
and
\begin{align*}
    \pi^{(i)}=\pi_i\circ\cdots\circ\pi_{1}:X_1\to X_{i+1}.
\end{align*}

For any $\mu\in \mathcal{M}^G(X_1)$, we write $\mu_i=\pi^{(i-1)}\mu\in \mathcal{M}^G(X_i)$.

Let $f: X_1\to \mathbb{R}$ be a continuous function and $d^{(i)}$ the metric on $X_i$. For any $F\in \mathcal{F}(G)$, the compatible metrics $d^{(i)}_F$ are defined as \eqref{el2.1} for all $i=1,\cdots,r$. Define a linear map from $C(X)$ to $\mathbb{R}$ by
\begin{equation*}
    \mathbb{S}_Ff(x)=\underset{g\in F}{\sum}f(gx).
\end{equation*}

Let \textbf{a}=$(a_1,\cdots,a_{r-1})\in \mathbb{R}^{r-1}$ with $0\leq a_i\leq 1$ for $i=1,\cdots,r-1$. Given $\Omega\subset X$ and $F\in \mathcal{F}(G)$, we define
\begin{flalign*}
& P^{\boldsymbol{a}}_1(\Omega, \spa f, \spa F, \spa \vep) &
\end{flalign*} \\[-35pt]
\begin{align*}
&= \inf \left\{ \spa \sum_{j=1}^n \exp \left( \sup_{U_j} S_F f \right) \spa \middle|
\begin{array}{l}
\text{$n \in \mathbb{N}$, $\{U_j\}_{j=1}^n$ is an open cover of $\Omega$} \\
\text{with $\diam(U_j, \spa d_F^{(1)}) < \vep$ for all $ 1 \spa \leq j \spa \leq n$}
\end{array}
\right\}.
\end{align*}
Suppose $P_i^{\textbf{a}}(\Omega, \spa f, \spa F, \spa \vep)$ is already defined, let $\Omega\subset X_{i+1}$, we define
\begin{flalign*}
& P^{\boldsymbol{a}}_{i+1}(\Omega, \spa f,\spa F, \spa \vep) &
\end{flalign*} \\[-35pt]
\begin{align*}
&= \inf \left\{ \spa \sum_{j=1}^n \Big( P^{\boldsymbol{a}}_i(\pi_i^{-1}(U_j), \spa f, \spa F, \spa \vep) \Big)^{a_i} \spa \middle|
\begin{array}{l}
\text{$n \in \mathbb{N}$, $\{U_j\}_{j=1}^n$ is an open cover of $\Omega$} \\
\text{with $\diam(U_j, \spa d_F^{(i+1)}) < \vep$ for all $ 1 \spa \leq j \spa \leq n$}
\end{array}
\right\}.
\end{align*}

\begin{definition}
    For each $\vep>0$, we define
    \begin{align*}
        {P}^{\textbf{a}}(X_r,f,\vep)={\lim_{n\to \infty}}\,\displaystyle\frac{\log\,{P}^{\textbf{a}}(X_r,f,F_n,\epsilon)}{|F_n|}.
    \end{align*}
\end{definition}
\begin{remark}

   For any $f\in C(X)$, $F\in \mathcal{F}(G)\mapsto \log\,{P}^{\textbf{a}}(X_r,f,F_n,\mathcal{U}_r)$ is $G$-invariant and sub-additive. Then by Theorem \ref{thm2.1} the limit
\begin{equation*}
    \underset{n\to \infty}{\lim}\,\displaystyle\frac{\log\,{P}^{\textbf{a}}(X_r,f,F_n,\vep)}{|F_n|}
\end{equation*}
exists and is independent of the choice of F$\o$lner sequence of $G$.
\end{remark}
Since ${P}^{\textbf{a}}(X_r,f,\vep)$ is non-decreasing as $\vep$ tends to $0$, we get the following definitions.
\begin{definition}
    The \textbf{topological pressure of a-exponent} $P^{\textbf{a}}(f,G)$ is defined by
    \[P^{\textbf{a}}(f,G)=\lim_{\epsilon\to 0}{P}^{\textbf{a}}(X_r,f,\vep).\]

If $f\equiv 0$, then the \textbf{topological entropy of a-exponent} $h^{\textbf{a}}_{top}(G)$ is given by
\[h^{\textbf{a}}_{top}(G)=P^{\textbf{a}}(0,G).\]
\end{definition}

Let \textbf{a}=$(a_1,a_2,\cdots,a_{r-1})\in \mathbb{R}^{r-1}$ with $0\leq a_i\leq 1$ for $i=1,\cdots,r-1$, we define a probability vector as in \cite{al}
$$ \boldsymbol{\omega}=\left\{
\begin{aligned}
w_1&=&a_1a_2a_3\cdots a_{r-1} \\
w_2&=&(1-a_1)a_2a_3\cdots a_{r-1} \\
\cdots\\
w_{r-1}&=&(1-a_{r-2})a_{r-1}\\
w_r&=&1-a_{r-1}
\end{aligned}
\right.
$$
Our main results are stated as follows.

\begin{theorem}\label{main}
   Let $(X_i,G)$ ($i=1,\cdots,r$) be dynamical systems, \textbf{a}=$(a_1,\cdots,a_{r-1})$ with $0\leq a_i\leq 1$ for $i=1,\cdots,r-1$ and $\pi_i:X_i\to X_{i+1}$. For any $f\in C(X_1)$, we have
    \begin{equation*}
        P^{\textbf{a}}(G,f)=\sup_{\mu\in\mathcal{M}^G(X_1)}\Big(\sum_{i=1}^rw_ih_{\mu_i}(X_i,G)+w_1\int_{X_1}fd\mu\Big),
    \end{equation*}
    where $\mu_i=\pi^{(i-1)}\mu$.
\end{theorem}

If $f\equiv 0$, we obtain the variational principle for weighted topological entropy, that is,
\begin{corollary}

    Let $(X_i,G)$ ($i=1,\cdots,r$) be dynamical systems, \textbf{a}=$(a_1,\cdots,a_{r-1})$ with $0\leq a_i\leq 1$ for all $i$ and $\pi_i:X_i\to X_{i+1}$, we have
    \begin{equation*}
        h^{\textbf{a}}(G)=\sup_{\mu\in\mathcal{M}^G(X_1)}\Big(\sum_{i=1}^rw_ih_{\mu_i}(X_i,G)\Big),
    \end{equation*}
    where $\mu_i=\pi^{(i-1)}\mu$.
\end{corollary}

For convenience, we write
\[P_{var}^{\textbf{a}}(G,f)=\sup_{\mu\in\mathcal{M}^G(X_1)}\Big(\sum_{i=1}^rw_ih_{\mu_i}(X_i,G)+w_1\int_{X_1}fd\mu\Big)\] and respectively,\[h^{\boldsymbol{a}}_{var}(G)=\sup_{\mu\in \mathcal{M}^G(X_1)}\Big(\sum_{i=1}^rw_ih_{\mu_i}(X_i,G)\Big)\].

\subsection{Measure-theoretical entropy}

Let $(X,\mathcal{B}_X,\mu)$ be a compact metrizable probability space where $\mu\in \mathcal{M}(X)$. Given two measurable finite paritions $\alpha,\beta\in \mathcal{P}_X$, the entropy of $\mathcal{A}$ is defined by
\[H_\mu(\alpha)=\sum_{A\in\alpha}-\mu(A)\log\mu(A),\]
and the conditional entropy of $\alpha$ with respect to $\beta$ is given by
\[H_\mu(\alpha|\beta)=-\sum_{B\in\beta}\mu(B)\sum_{A\in \alpha}\displaystyle\frac{\mu(A\cap B)}{\mu(B)}\log\displaystyle\frac{\mu(A\cap B)}{\mu(B)}.\]

To prove the variational principle, unlike $\mathbb{Z}$-action in \cite{mt,al}, it seems hard to use an amplification trick for amenable groups. Hence, our proof needs to borrow some techniques of the local entropy theory (for example, see \cite{ppr,hyz,ly}). First, we recall some notations and properties below.

Let $\mu\in \mathcal{M}^G(X)$, given a measurable cover $\mathcal{U}\in \mathcal{C}_X$, we define the $\mu$-entropy of $\mathcal{U}$ as
\[H_\mu(\mathcal{U})=\inf_{\alpha\in \mathcal{P}_X}H_\mu(\alpha).\]
For any $\alpha\in \mathcal{P}_X$, we define
\[h_{\mu}(\alpha,G)=\lim_{n\to \infty}\frac{H_\mu(\alpha_{F_n})}{|F_n|}=\inf_{F\in \mathcal{F}(G)}\frac{H_{\mu}(\alpha_F)}{|F|},\]
 see \cite[Lemma 3.1]{hyz} for the second equality.
For any $\mathcal{U}\in \mathcal{C}^o_X$, we define
\[h_\mu^+(\mathcal{U},G)=\inf_{\mathcal{U}\preceq\alpha\in\mathcal{P}_X}\Big(\lim_{n\to \infty}\displaystyle\frac{H_\mu(\alpha_{F_n})}{|F_n|}\Big) \text{ and } h_\mu^-(\mathcal{U},G)=\lim_{n\to \infty}\displaystyle\frac{H_\mu(\mathcal{U}_{F_n})}{|F_n|}.\]
It is easy to see that $h_\mu^+(\mathcal{U},G)\geq h_\mu^-(\mathcal{U},G)$. Moreover, the authors in \cite{hyz} proved the following significant result:
\begin{lemma}\rm{(\cite[Theorem 4.14]{hyz})}\label{lem2.7}
   Let $\mu\in \mathcal{M}^G(X)$, for any $\mathcal{U}\in\mathcal{C}_X^o$, we have
   \[h_\mu^+(\mathcal{U},G)=h_\mu^-(\mathcal{U},G).\]

\end{lemma}

To prove the variational principle, we need a Lemma in \cite{hyz}.
\begin{lemma}\rm{(\cite[Lemma 3.1 (3)]{hyz}\label{lem2.13})} Let $\alpha\in \mathcal{P}_X$, $\mu\in \mathcal{M}^G(X)$ and $A,F\in \mathcal{F}(G)$, then we have
    \[H_\mu(\alpha_F)\leq\sum_{g\in F}\frac{1}{|A|}H_\mu(\alpha_{Ag})+|F\backslash\{g\in G:A^{-1}g\subset F\}|\cdot\log |\alpha|,\]
    where $|\cdot|$ denotes the cardinal operation. Hence, if $e_G\in A$
     \[H_\mu(\alpha_F)\leq\sum_{g\in F}\frac{1}{|A|}H_\mu(\alpha_{Ag})+|B(F,A^{-1})|\cdot\log |\alpha|,\]
\end{lemma}

\begin{definition}
The measure-theoretical entropy of $(X,G)$ is defined by
    \begin{equation}\label{el2.2}
    h_\mu(X,G)=\sup_{\alpha\in \mathcal{P}_X}h_\mu(\alpha,G)=\sup_{\mathcal{U}\in\mathcal{C}^o_X}h_\mu^+(\mathcal{U},G)=\sup_{\mathcal{U}\in\mathcal{C}^o_X}h_\mu^-(\mathcal{U},G).
\end{equation}
\end{definition}

\begin{remark}\label{rmk2.9}
\begin{enumerate}
    \item  The first equality in \eqref{el2.2} is the Kolmogorov-Sinai entropy for amenable group actions;
    \item Moreover, from \cite[Theorem 3.5]{hyz}, we have
    $h_\mu(X,G)=\sup_{\mathcal{U}\in \mathcal{C}^o_X}h^{\pm}(\mathcal{U},G).$
    Therefore, the second equality is obtained.
\end{enumerate}
\end{remark}

\subsection{Zero-dimensional principal extensions} \label{subsection: zero-dimensional principal extension}
The subsection aims to review the definitions and properties of zero-dimensional principal extensions, one can refer to Tsukamoto's paper \cite{mt} and the book of Downarowicz \cite{Downarowicz}. This part will be useful to show that if the given dynamical systems are zero-dimensional, $P^{\boldsymbol{a}}(G,f) \leq P^{\boldsymbol{a}}_\mathrm{var}(G,f)$ in Section 4.

Before we discuss the principle extensions, we need some notations.
Let $(X_i, \spa G)$ $(i=1, \spa 2, \spa \ldots, \spa r)$ be $G$-systems, $\pi_i: X_i \rightarrow X_{i+1} \hspace{5pt} (i=1, \spa 2, \spa \ldots \spa , \spa r-1)$ being factor maps, $\boldsymbol{a} = (a_1, \cdots, a_{r-1}) \in [0, 1]^{r-1}$, and $f: X_1 \to \mathbb{R}$ a continuous function.

\begin{definition}\rm{(\cite[Definition 3.1]{al})} \label{definition: chain of covers}
Consider a cover $\flo{U}^{(i)}$ of $X_i$ for each $i$. For an F$\o$lner sequence $\{F_n\}$ and a positive number $\vep$, the family $(\flo{U}^{(i)})_i$ is said to be \textbf{a chain of ($\boldsymbol{F_n}$, $\boldsymbol{\vep}$)-covers} of $(X_i)_i$ if the following conditions are satisfied:
\begin{enumerate}
\item[(1)] For every $i$ and $V \in \flo{U}^{(i)}$, we have $\diam(V, d^{(i)}_{F_n}) < \vep$.
\item[(2)] For each $1 \leq i \leq r-1$ and $U \in \flo{U}^{(i+1)}$, there is $\flo{U}^{(i)}(U) \subset \flo{U}^{(i)}$ such that
\[ \pi_i^{-1}(U) \subset \bigcup \flo{U}^{(i)}(U) \]
and
\[ \flo{U}^{(i)} = \bigcup_{U \in \flo{U}^{(i+1)}} \flo{U}^{(i)}(U). \]
\end{enumerate}
Moreover, if all the elements of each $\flo{U}^{(i)}$ are open/closed/compact, we call $(\flo{U}^{(i)})_i$ \textbf{a chain of open/closed/compact ($\boldsymbol{F_n}$, $\boldsymbol{\vep}$)-covers} of $(X_i)_i$.
\end{definition}

\begin{remark} \label{remark: chains of covers}
Note that we can rewrite $P^{\boldsymbol{a}}_r(X_r, \spa f, \spa F_n, \spa \vep)$ using chains of open covers as follows. For a chain of ($F_n$, $\vep$)-covers $(\flo{U}^{(i)})_i$ of $(X_i)_i$, let \\[-10pt]
\begin{flalign*}
& \flo{P}^{\boldsymbol{a}}\left( f, \spa F_n, \spa \vep, \spa (\flo{U}^{(i)})_i \right) &
\end{flalign*} \\[-35pt]
\begin{align*}
&= \sum_{U^{(r)} \in \flo{U}^{(r)}} \left( \sum_{U^{(r-1)} \in \flo{U}^{(r-1)}(U^{(r)})} \left( \cdots \left( \sum_{U^{(1)} \in \flo{U}^{(1)}(U^{(2)})} e^{\sup_{U^{(1)}}S_{F_n}f} \right)^{a_1} \cdots \right)^{a_{r-2}}\right)^{a_{r-1}}.
\end{align*}
Then
\begin{flalign*}
& P^{\boldsymbol{a}}_r(X_r, \spa f, \spa F_n, \spa \vep) &
\end{flalign*}
\begin{align*}
&= \inf{ \left\{ \flo{P}^{\boldsymbol{a}}\left( f, \spa F_n, \spa \vep, \spa (\flo{U}^{(i)})_i \right) \spa \middle| \spa (\flo{U}^{(i)})_i \text{ is a chain of open ($F_n$, $\vep$)-covers of $(X_i)_i$ } \right\} }.
\end{align*}
\end{remark}

Suppose $\pi: (Y, \spa G) \rightarrow (X, \spa G)$ is a factor map between $G$-systems. Let $d$ be a metric on $Y$. We define the \textbf{conditonal topological entropy} of $\pi$ by
\[ h_\mathrm{top}(G,\spa Y| \spa X) = \lim_{\vep \to 0} \left( \lim_{n \to \infty} \frac{\sup_{x \in X} \log{\#(\pi^{-1}(x), F_n, \vep)}}{|F_n|} \right). \]
Here,
\begin{align*}
\#(\pi^{-1}(x), \spa F_n, \spa \vep)
&= \min \left\{ \spa m \in \mathbb{N} \spa \middle|
\begin{array}{l}
\text{There exists an open cover $\{U_j\}_{j=1}^m$ of $\pi^{-1}(x)$} \\
\text{such that $\diam(U_j, \spa d_{F_n}) < \vep$ for all $ 1 \spa \leq j \spa \leq m$}
\end{array}
\right\}.
\end{align*}
it is easy to check that $F\in\mathcal{F}(G)\mapsto\sup_{x \in X} \log{\#(\pi^{-1}(x), F, \vep)}$ is sub-additive, one can see \cite{yk} for more details about conditional topological entropy for amenable group actions.

A factor map $\pi: (Y, \spa G) \rightarrow (X, \spa G)$ between $G$-systems is said to be a \textbf{principal factor map} if
\[ h_\mathrm{top}(G,\spa Y| \spa X) = 0. \]
Also, $(Y, \spa G)$ is called a \textbf{principal extension} of $(X, \spa G)$.

The following theorem is from \cite[Theorem 2.5]{YCYY}. It says that a factor map is principal if and only if it preserves measure-theoretic entropy.
\begin{theorem}
Suppose $\pi: (Y, \spa G) \rightarrow (X, \spa G)$ is a principal factor map. Then $\pi$ preserves measure-theoretic entropy, namely,
\[ h_\mu(Y,G) = h_{\pi*\mu}(X,G) \]
for any $G$-invariant probability measure $\mu$ on Y.
\end{theorem}

Suppose $\pi: (X_1, G) \rightarrow (X_2, G)$ and $\phi: (Y, G) \rightarrow (X_2, G)$ are factor maps between $G$-systems. A fiber product $(X_1 \times_{X_2}^{} Y, \spa G)$ of $(X_1, G)$ and $(Y, G)$ over $(X_2, G)$ is defined by
\[ X_1 \times_{X_2}^{} Y = \left\{ (x, y) \in X_1 \times Y \spa \middle| \spa \pi(x) = \phi(y) \right\}, \]
\[ g: X_1 \times_{X_2} Y \ni (x, y) \longmapsto \left(gx, gy\right) \in X_1 \times_{X_2}^{} Y,\text{ for every $g\in G$. }\]
From the definition, one can get the following commutative diagram:
\begin{equation} \label{diagram: row extension one}
\begin{gathered}
 \xymatrix@C=40pt@R=36pt{
	X_1 \times_{X_2}^{} Y \ar[d]_-{\pi'} \ar[r]^-\psi & X_1 \ar[d]^-\pi \\
	Y \ar[r]_-\phi & X_2 },
\end{gathered}
\end{equation}
here $\pi'$ and $\psi$ are restrictions of the projections onto $Y$ and $X_1$, respectively:
\[ \pi':  X_1 \times_{X_2}^{} Y \ni (x, y) \longmapsto y \in Y, \]
\[ \psi: X_1 \times_{X_2}^{} Y \ni (x, y) \longmapsto x \in X_1. \]
Since $\pi$ and $\phi$ are surjective, both $\pi'$ and $\psi$ are factor maps.

\begin{lemma} \label{lemma: square extension}
If $\phi$ is a principal extension in the diagram (\refeq{diagram: row extension one}), then $\psi$ is also a principal extension.
\end{lemma}
\begin{proof}
The proof of the result is similar to the process of \cite[Lemma 3.7]{al} and \cite[Lemma 5.3]{mt}.
\end{proof}

A dynamical system $(Y, G)$ is said to be \textbf{zero-dimensional} if it is totally disconnected, equivalent to say that the open basis of the topology of $Y$ is formed by clopen subsets. For example, The Cantor set $\{ 0, 1 \}^\mathbb{N}$ with the shift map is such a zero-dimensional dynamical system.

A principal extension $(Y, \spa G)$ of $(X, \spa G)$ is called a \textbf{zero-dimensional principal extension} if $(Y, G)$ is zero-dimensional. The following important theorem can be found in \cite[Theorem 3.2]{dh1}.
\begin{theorem} \label{theorem: zero dimensional principal extension}
For any $G$-system, there is a zero-dimensional principal extension.
\end{theorem}

Let $(Y_i, \spa G)$ ($i=1, \spa 2, \spa \ldots, \spa m$) be $G$-systems, $\pi_i: Y_i \rightarrow Y_{i+1} \hspace{5pt} (i=1, \spa 2, \spa \cdots \spa , \spa m-1)$ factor maps, and $\boldsymbol{a} = (a_1, \cdots, a_{m-1}) \in [0, 1]^{m-1}$. Fix $2 \leq k \leq m-1$, we choose a zero-dimensional principal extension $\phi_k: (Z_k, G) \rightarrow (Y_k, G)$.
For each $1 \leq i \leq k-1$, let $(Y_i \times_{Y_k} Z_k, G)$ be the fiber product and $\phi_i: Y_i \times_{Y_k} Z_k \rightarrow Y_i$ be the restriction of the projection as in the earlier definition. Then we have
\[ \xymatrix@C=60pt{
Y_i \times_{Y_k} Z_k \ar[r]^-{\phi_i} \ar[d] & Y_i \ar[d]^-{\pi_{k-1} \circ \pi_{k-2} \circ \cdots \circ \pi_i} \\
Z_k \ar[r]_-{\phi_k} & Y_k } \]
By Lemma \ref{lemma: square extension}, $\phi_i$ is a principal factor map.
We define $\Pi_i: Y_i \times_{Y_k} Z_k \rightarrow Y_{i+1} \times_{Y_k} Z_k$ by $\Pi_i(x, y) = \left( \pi_i(x), y \right)$ for each $i$. Then we have the following commutative diagram:
\begin{equation} \label{diagram: line extension}
\begin{gathered}
\xymatrix@C=18pt@R=60pt{
Y_1 \times_{Y_k}^{} Z_k \ar[r]^-{\Pi_1} \ar[d]_-{\phi_1} & Y_2 \times_{Y_k}^{} Z_k\ar[r]^-{\Pi_2} \ar[d]^-{\phi_2}&\cdots \ar[r]^-{\Pi_{k-2}}&Y_{k-1} \times_{Y_k}^{} Z_k\ar[r]^-{\Pi_{k-1}} \ar[d]_-{\phi_{k-1}}&Z_k\ar[d]_-{\phi_k}\ar[rd]^-{\pi_k \circ \phi_k}\\
Y_1\ar[r]^-{\pi_1}&Y_2 \ar[r]^-{\pi_2}&\cdots\ar[r]^-{\pi_{k-2}}& Y_{k-1} \ar[r]^-{\pi_{k-1}}& Y_k \ar[r]^-{\pi_k}& Y_{k+1} \ar[r]^-{\pi_{k+1}}& \cdots \ar[r]^-{\pi_{m-1}}& Y_m}.
\end{gathered}
\end{equation}
Let
\begin{equation*}
\begin{gathered}
(Z_i, G) = (Y_i \times^{}_{Y_k} Z_k, G) \text{ for $1 \leq i \leq k-1$},
\hspace{6pt} (Z_i, G) = (Y_i, G) \text{   for $k+1 \leq i \leq m$}, \\
\Pi_k = \pi_k \circ \phi_k: Z_k \rightarrow Y_{k+1}, \hspace{6pt}
\Pi_i = \pi_i: Z_i \rightarrow Z_{i+1} \text{  for $k+1 \leq i \leq m-1$}, \\
\phi_i = \mathrm{id}_{Z_i}: Z_i \rightarrow Z_i  \text{   for $k+1 \leq i \leq m$}.
\end{gathered}
\end{equation*}

\begin{lemma} \label{lemma: pressure inequality}
In the settings above,
\begin{equation}\label{eq2.1}
 P^{\boldsymbol{a}}_{\mathrm{var}}(f, G, \boldsymbol{\pi}) \geq P^{\boldsymbol{a}}_{\mathrm{var}}(f \circ \phi_1, G, \boldsymbol{\Pi})
\end{equation}
and
\begin{equation}\label{eq2.2}
P^{\boldsymbol{a}}(f, G, \boldsymbol{\pi}) \leq P^{\boldsymbol{a}}(f \circ \phi_1, G, \boldsymbol{\Pi}).
\end{equation}
Here, $\boldsymbol{\pi} = (\pi_i)_i$, $\boldsymbol{\Pi} = (\Pi_i)_i$.
\end{lemma}
\begin{proof}
We note that the following process does not need $Z_k$ to be zero-dimensional. Let
\begin{gather*}
\pi^{(0)} = \mathrm{id}_{Y_1}: Y_1 \to Y_1, \\
\pi^{(i)} = \pi_i \circ \pi_{i-1} \circ \cdots \circ \pi_1: Y_1 \to Y_{i+1}
\end{gather*}
and
\begin{gather*}
\Pi^{(0)} = \mathrm{id}_{Z_1}: Z_1 \to Z_1, \\
\Pi^{(i)} = \Pi_i \circ \Pi_{i-1} \circ \cdots \circ \Pi_1: Z_1 \to Z_{i+1}.
\end{gather*}
Let $\nu \in \mathcal{M}^G(Z_1)$ and $1 \leq i \leq m$. Since all the vertical maps in \eqref{diagram: line extension} are principal factor maps, we have
\[ h_{{\Pi^{(i-1)}}*\nu}(Z_i,G) = h_{({\phi_i}\circ{\Pi^{(i-1)}})*\nu}(Y_i,G) = h_{({\pi^{(i-1)}}\circ\phi_1)*\nu}(Y_i,G). \]
Then we get that
\begin{align*}
P^{\boldsymbol{a}}_{\mathrm{var}}(f \circ \phi_1,G, \boldsymbol{\Pi})
&= \sup_{\nu \in \mathcal{M}^G(Z_1)} \left( \sum_{i=1}^{m} w_i h_{{\Pi^{(i-1)}}*\nu} (Z_i,G) + w_1 \int_{Z_1} f \circ \phi_1 d\nu \right) \\
&= \sup_{\nu \in \mathcal{M}^G(Z_1)} \left( \sum_{i=1}^{m} w_i h_{({\pi^{(i-1)}}\circ\phi_1)*\nu}(Y_i,G) + w_1 \int_{Y_1} f d\big( (\phi_1)*\nu \big) \right) \\
&\leq \sup_{\mu \in \mathcal{M}^G(Y_1)} \left( \sum_{i=1}^{m} w_i h_{{\pi^{(i-1)}}*\mu} (Y_i,G) + w_1 \int_{Y_1} f d\mu \right) \\[4pt]
&= P^{\boldsymbol{a}}_{\mathrm{var}}(f, G, \boldsymbol{\pi}).
\end{align*}
Then the first equation \eqref{eq2.1} is proved.

Let $d^i$ be a metric on $Y_i$ for each $i$ and $\widetilde{d^k}$ a metric on $Z_k$. We define a metric $\widetilde{d^i}$ on $(Z_i, G)$ for $1 \leq i \leq k-1$ by
\[ \widetilde{d^i}\big( (x_1, y_1), (x_2, y_2) \big) = \max{\{ d^i(x_1, x_2), \widetilde{d^k}(y_1, y_2) \}}.\]
where $(x_1, y_1), (x_2, y_2) \in Z_i = Y_i \times_{Y_k} Z_k$ and we set $\widetilde{d^i} = d^i$ for $k+1 \leq i \leq m$.

Let $\vep>0$. By the continuity of $\phi$, we can choose $0 < \delta < \vep$ such that for every $1 \leq i \leq m$, if $\widetilde{d^i}(x, y) < \delta$ with $x, y \in Z_i$,  $d^i( \phi_i(x), \phi_i(y) ) < \vep$.

Let $\{F_n\}$ be a F$\o$lner sequence. We claim that
\[ P^{\boldsymbol{a}}_r(f, G, \boldsymbol{\pi}, F_n, \vep)
\leq P^{\boldsymbol{a}}_r(f \circ \phi_1, G, \boldsymbol{\Pi}, F_n, \delta). \]
Take $M > 0$ with
\[
P^{\boldsymbol{a}}_r(f \circ \phi_1, G, \boldsymbol{\Pi}, F_n, \delta) < M.
\]
Then there exists a chain of open ($F_n$, $\delta$)-covers $(\flo{U}^{(i)})_i$ of $(Z_i)_i$ (see Definition \ref{definition: chain of covers} and Remark \ref{remark: chains of covers}) with
\[ \flo{P}^{\boldsymbol{a}}\left(f \circ \phi_1, \spa G, \spa \boldsymbol{\Pi}, \spa F_n, \spa \delta, \spa (\flo{U}^{(i)})_i \right) < M. \]
For each $U \in \flo{U}^{(m)}$, We can find a compact set $C_U \subset U$  such that $\bigcup_{U \in \flo{U}^{(m)}} C_U = Z_m$. Let $\flo{K}^{(m)} := \{ C_U \spa | \spa U \in \flo{U}^{(m)} \}$. Since $\Pi_{m-1}^{-1}(C_U) \subset \Pi_{m-1}^{-1}(U)$ is compact for each $U \in \flo{F}^{(m)}$, we can find a compact set $E_V \subset V$ for each $V \in \flo{U}^{(m-1)}(U)$ such that \[
\Pi_{m-1}^{-1}(C_U) \subset \bigcup_{V \in \flo{U}^{(k)}(U)} E_V.
\]
Let
\[\flo{K}^{(m-1)}(C_U) := \{ E_V \spa | \spa V \in \flo{U}^{(m-1)}(U) \}\]
and
\[
\flo{K}^{(m-1)} := \bigcup_{C \in \flo{K}^{(m)}} \flo{K}^{(m-1)}(C).
\]
Continue the process and we can obtain a chain of compact ($F_n$, $\delta$)-covers $(\flo{K}^{(i)})_i$ of $(Z_i)_i$ with
\[ \flo{P}^{\boldsymbol{a}}\left(f \circ \phi_1, \spa G, \spa \boldsymbol{\Pi}, \spa F_n, \spa \delta, \spa (\flo{K}^{(i)})_i \right)
\leq \flo{P}^{\boldsymbol{a}}\left(f \circ \phi_1, \spa G, \spa \boldsymbol{\Pi}, \spa F_n, \spa \delta, \spa (\flo{U}^{(i)})_i \right) < M. \]
Let $\phi_i(\flo{K}^{(i)}) = \left\{ \phi_i(C) \spa \middle| \spa C \in \flo{K}^{(i)} \right\}$ for each $i$. Note that for any $\Omega \subset Z_i$,
\[ \pi_{i-1}^{-1} ( \phi_{i} ( \Omega ) ) = \phi_{i-1} ( \Pi_{i-1}^{-1} ( \Omega ) ). \]
Therefore, $(\phi_i(\flo{K}^{(i)}))_i$ is a chain of compact ($F_n$, $\vep$)-covers of $(Y_i)_i$. We have
\begin{align*}
\flo{P}^{\boldsymbol{a}}\left( f, \spa G, \spa \boldsymbol{\pi}, \spa F_n, \spa \vep, \spa (\phi_i(\flo{K}^{(i)}))_i \right)
&\leq
\flo{P}^{\boldsymbol{a}}\left( f \circ \phi_1, \spa G, \spa \boldsymbol{\Pi}, \spa F_n, \spa \delta, \spa (\flo{K}^{(i)})_i \right) < M.
\end{align*}

Since $f$ is continuous and each $\phi_i(\flo{K}^{(i)})$ is a closed cover, we can choose a chain of open ($F_n$, $\vep$)-covers $(\flo{O}^{(i)})_i$ of $(Y_i)_i$ satisfying that for every $i$, $\phi_i(\flo{K}^{(i)})\subseteq \flo{O}^{(i)}$ and
\[
\flo{P}^{\boldsymbol{a}}\left( f, \spa G, \spa \boldsymbol{\pi}, \spa F_n, \spa \vep, \spa (\flo{O}^{(i)})_i \right) < M.
\]
Therefore,
\[
P^{\boldsymbol{a}}_r(f, G, \boldsymbol{\pi}, F_n, \vep) \leq \flo{P}^{\boldsymbol{a}}\left( f, \spa G, \spa \boldsymbol{\pi}, \spa F_n, \spa \vep, \spa (\flo{O}^{(i)})_i \right) < M.
\]
Since $M$ was chosen arbitrarily, we have
\[
P^{\boldsymbol{a}}_r(f, G, \boldsymbol{\pi}, F_n, \vep) \leq P^{\boldsymbol{a}}_r(f \circ \phi_1, G, \boldsymbol{\Pi}, F_n, \delta).
\]
This implies
 \[
 P^{\boldsymbol{a}}(f, G, \boldsymbol{\pi}) \leq P^{\boldsymbol{a}}(f \circ \phi_1, G, \boldsymbol{\Pi}).
 \]
\end{proof}

The following proposition play a central role in the proof of $P^{\boldsymbol{a}}(f,G) \leq P^{\boldsymbol{a}}_\mathrm{var}(f,G)$ in Section 4 when all dynamical systems are zero-dimensional.

\begin{proposition} \label{proposition: zero-dimensional trick}
For all $G$-systems $(X_i, \spa G)$ ($i=1, \spa 2, \spa \ldots, \spa r$) and factor maps $\pi_i: X_i \rightarrow X_{i+1} \hspace{5pt} (i=1, \spa 2, \spa \cdots \spa , \spa r-1)$, there are zero-dimensional $G$-systems $(Z_i, \spa G)$ ($i=1, \spa 2, \spa \ldots, \spa r$) and factor maps $\Pi_i: Z_i \rightarrow Z_{i+1} \hspace{5pt} (i=1, \spa 2, \spa \cdots \spa , \spa r-1)$ such that the following property is satisfied: for every continuous function $f: X_1 \rightarrow \mathbb{R}$, there exists a continuous function $g: Z_1 \rightarrow \mathbb{R}$ with
\[ P^{\boldsymbol{a}}_{\mathrm{var}}(f, G, \boldsymbol{\pi}) \geq P^{\boldsymbol{a}}_{\mathrm{var}}(g, G, \boldsymbol{\Pi}) \]
and
\[ P^{\boldsymbol{a}}(f, G, \boldsymbol{\pi}) \leq P^{\boldsymbol{a}}(g, G, \boldsymbol{\Pi}). \]
\end{proposition}

\begin{proof}
\textbf{Step 1:} Construct zero-dimensional $G$-systems $(Z_i, \spa G)$ ($i=1, \spa 2, \spa \ldots, \spa r$) and factor maps $\Pi_i: Z_i \rightarrow Z_{i+1} \hspace{5pt} (i=1, \spa 2, \spa \cdots \spa , \spa r-1)$ alongside the following commutative diagram of $G$-systems and factor maps:

\begin{equation} \label{diagram: Xi to Zi}
\xymatrix@C=9pt@R=6pt{
	Z_1 \ar[r]^-{\phi_r} \ar[rdd]_{\Pi_1} & \ar[dd] & & & \cdots \hspace{90pt} \ar[r]^-{\phi_2} & X_1 \times_{X_r}^{} Z_r \ar[r]^-{\phi_1} \ar[dd]^{\pi_1^{(2)}} & X_1 \ar[dd]^{\pi_1} \\
	&&&&&& \\
	      & Z_2 \ar[r] \ar[rdd]_{\Pi_2} & & & \cdots \hspace{90pt} \ar[r]  & X_2 \times_{X_r}^{} Z_r \ar[r] \ar[dd]^{\pi_2^{(2)}} & X_2 \ar[dd]^{\pi_2} \\
	&&&&&& \\
	& & \ddots \hspace{4pt} \ar[rd]_{\Pi_{r-3}} & \ar[d]^{\pi^{(4)}_{r-3}} & \vdots \ar[d]^{\pi_{r-3}^{(3)}} & \vdots \ar[d]^{\pi_{r-3}^{(2)}} & \vdots \ar[d]^{\pi_{r-3}} \\
	      & & &	Z_{r-2} \ar[r] \ar[rddd]_{\Pi_{r-2}} & \left(X_{r-2} \times_{X_r}^{} Z_r\right) \times_{(X_{r-1} \times_{X_r} Z_r)}^{} Z_{r-1} \ar[r] \ar[ddd]^{\pi_{r-2}^{(3)}} & X_{r-2} \times_{X_r}^{} Z_r \ar[r] \ar[ddd]^{\pi_{r-2}^{(2)}} & X_{r-2} \ar[ddd]^{\pi_{r-2}} \\
	&&&&&& \\ 	&&&&&& \\
	      & & & & 	    Z_{r-1} \ar[r]^-{\psi_{r-1}} \ar[rdddd]_{\Pi_{r-1}} & X_{r-1} \times_{X_r}^{} Z_r \ar[r] \ar[dddd]^{\pi_{r-1}^{(2)}} & X_{r-1} \ar[dddd]^{\pi_{r-1}} \\
	&&&&&& \\	&&&&&& \\ 	&&&&&& \\
	& & & & & Z_r \ar[r]^-{\psi_r} \ar[rdd] & \hspace{2pt} X_r \ar[dd] \\
	&&&&&& \\
	& & & & & & \{*\} }
\end{equation}
where all the horizontal maps are principal factor maps.

By Theorem \ref{theorem: zero dimensional principal extension}, there is a zero-dimensional principal extension $\psi_r: (Z_r, G) \rightarrow (X_r, G)$. Denote the $G$-system consists of one point (the trivial $G$-system) by $\{*\}$. one can find that the maps $X_r \rightarrow \{*\}$ and $Z_r \rightarrow \{*\}$ send every element to $*$. For each $1 \leq i \leq r-1$, by Lemma \ref{lemma: square extension}, we can get that the map $X_i \times_{X_r}^{} Z_r \rightarrow X_i$ in the following diagram
\[ \xymatrix@C=40pt@R=36pt{
	X_i \times_{X_r}^{} Z_r \ar[d] \ar[r] & X_i \ar[d]^-{\pi_{r-1} \circ \pi_{r-2} \circ \cdots \circ \pi_i} \\
	Z_r \ar[r]_-{\psi_r} & X_r } \]
is a principal factor map.

For $1 \leq i \leq r-2$, define $\pi_i^{(2)}: X_i \times_{X_r}^{} Z_r \rightarrow X_{i+1} \times_{X_r}^{} Z_r$ by $\pi_i^{(2)}(x, z) = (\pi_i(x), y).$
Then every horizontal map in the right two rows of (\refeq{diagram: Xi to Zi}) is a principal factor map.

\textbf{Step 2:} Set a zero-dimensional principal extension $\psi_{r-1}: (Z_{r-1}, G) \rightarrow (X_{r-1} \times_{X_r}^{} Z_r, G)$ and let $\Pi_{r-1} = \pi_{r-1}^{(2)} \circ \psi_{r-1}$. The rest of (\refeq{diagram: Xi to Zi}) is constructed similarly, and by Lemma \ref{lemma: square extension}, each horizontal map is a principal factor map.

Let $f: X_1 \rightarrow \mathbb{R}$ be a continuous map. Applying Lemma \ref{lemma: pressure inequality} to the right two rows of (\refeq{diagram: Xi to Zi}), we get
\[ P^{\boldsymbol{a}}_{\mathrm{var}}(f, G, \boldsymbol{\pi}) \geq P^{\boldsymbol{a}}_{\mathrm{var}}(f \circ \phi_1, G, \boldsymbol{\Pi^{(2)}}) \]
and
\[ P^{\boldsymbol{a}}(f, G, \boldsymbol{\pi}) \leq P^{\boldsymbol{a}}(f \circ \phi_1, G, \boldsymbol{\Pi^{(2)}}) \]
for $\boldsymbol{\Pi^{(2)}} = (\pi^{(2)}_i)_i$. Again by Lemma \ref{lemma: pressure inequality},
\[ P^{\boldsymbol{a}}_{\mathrm{var}}(f \circ \phi_1, G, \boldsymbol{\Pi^{(2)}}) \geq P^{\boldsymbol{a}}_{\mathrm{var}}(f \circ \phi_1 \circ \phi_2, G, \boldsymbol{\Pi^{(3)}}) \]
and
\[ P^{\boldsymbol{a}}(f \circ \phi_1, G, \boldsymbol{\Pi^{(2)}}) \leq P^{\boldsymbol{a}}(f \circ \phi_1 \circ \phi_2, G, \boldsymbol{\Pi^{(3)}}) \]
where $\boldsymbol{\Pi^{(3)}} = \big( (\pi^{(3)}_i)_{i=1}^{r-2}, \Pi_{r-1} \big)$.
We continue inductively and obtain the desired inequalities, where $g$ is taken as $f \circ \phi_1 \circ \phi_2 \circ \cdots \circ \phi_r$.
\end{proof}
\section{Proof of $P_{var}^{\textbf{a}}(G,f)\leq P^{\textbf{a}}(G,f)$ and $h_{var}^{\textbf{a}}(G)\leq h_{top}^{\textbf{a}}(G)$}
The following theorem gives a half variation principle for the higher dimensional topological pressure.
\begin{theorem}\label{thm3.1}
    Let $(X_i, G)$, $i=1,\cdots,r$, be dynamical systems with $G$ being a discrete countable amenable group and $\pi_{i}: X_i\to X_{i+1}$, $i=1,\cdots,r-1$ factor maps. Then we obtain the \textit{lower bound}
    \[P_{var}^{\textbf{a}}(G,f)\leq P^{\textbf{a}}(G,f),\]
    for any $f\in C(X)$.
\end{theorem}
When $f\equiv0$, as a direct consequence of Theorem \ref{thm3.1}, we have the following result.
\begin{corollary}\label{cor3.2}
     Let $(X_i, G)$, $i=1,\cdots,r$, be dynamical systems with $G$ a discrete countable amenable group and $\pi_{i}: X_i\to X_{i+1}$, $i=1,\cdots,r-1$ factor maps. Then we obtain the \textit{lower bound}
     \[h_{var}^{\textbf{a}}(G)\leq h^{\textbf{a}}(G).\]
\end{corollary}
Before proving Theorem \ref{thm3.1}, We need a following classical Lemma.
\begin{lemma}\rm{(\cite[Lemma 9.9]{pw}\label{lem3.3})}
    Let $a_1,\cdots,a_n$ be given real numbers. If $p_i\geq 0$, $i=1,\cdots,n$, and $\sum_{i=1}^np_i=1$, then
    \[\sum_{i=1}^n(-p_i\log p_i+p_ia_i)\le\log\sum_{i=1}^n e^{a_i}.\]
\end{lemma}

Now we start to prove the half principle.
\begin{proof}[Proof of Theorem \ref{thm3.1}]
To make our proof clearer, we need to define some new quantities as follows.

Let $\mathcal{U}^{(i)}\in \mathcal{C}^o_X$ be open covers of $X_i$, respectively, for any $\Omega_1\subset X_1$, we define
\begin{flalign*}
& P^{\boldsymbol{a}}_1(\Omega_1, \spa f, \spa F, \spa \mathcal{U}^{(1)}) &
\end{flalign*} \\[-35pt]
\begin{align}
&= \inf \left\{ \spa \sum_{j=1}^n \exp \left( \sup_{U_j} S_F f \right) \spa \middle|
\begin{array}{l}
\text{$n \in \mathbb{N}$, $\{V_j\}_{j=1}^n$ is a cover of $\Omega$} \\
\text{with $V_i\in\mathcal{V}\succeq\mathcal{U}^{(1)}_F$  for some $\mathcal{V}\in \mathcal{C}_{X_1}$ }
\end{array}
\right\}.
\end{align}\label{e3.1}
Suppose $P_i^{\textbf{a}}(\Omega_i, \spa f, \spa F, \spa \mathcal{U}^{(i)})$ is already defined, let $\Omega_{i+1}\subset X_{i+1}$, we define
\begin{flalign*}
& P^{\boldsymbol{a}}_{i+1}(\Omega_{i+1}, \spa f,\spa F, \spa \mathcal{U}^{(i+1)}) &
\end{flalign*} \\[-35pt]
\begin{align}\label{e3.2}
&= \inf \left\{ \spa \sum_{j=1}^n \Big( P^{\boldsymbol{a}}_i(\pi_i^{-1}(V_j), \spa f, \spa F, \spa \mathcal{U}^{(i)}) \Big)^{a_i} \spa \middle|
\begin{array}{l}
\text{$n \in \mathbb{N}$, $\{V_j\}_{j=1}^n$ is a cover of $\Omega_{i+1}$} \\
\text{with $V_j\in \mathcal{V}\succeq\mathcal{U}^{(i+1)}_F$ for some $\mathcal{V}\in \mathcal{C}_{X_{i+1}}$}
\end{array}
\right\}.
\end{align}
\begin{remark}\label{remark3.4}
    Note that the cover $\mathcal{V}$ in \eqref{e3.2} can be chosen from partitions for all $i=1,\cdots,r$. If $\{V_1,\cdots, V_n\}$ is a cover of $\Omega_{i+1}$, then $\{W_i: W_i=V_i\backslash\cup_{j<i}V_j,i=1,\cdots,n\}$ is also a cover of $\Omega_{i+1}$ such that
\[\sum_{j=1}^n \Big( P^{\boldsymbol{a}}_i(\pi_i^{-1}(V_j), \spa f, \spa F, \spa \mathcal{U}^{(i)}) \Big)^{a_i}\geq \sum_{j=1}^n \Big( P^{\boldsymbol{a}}_i(\pi_i^{-1}(W_j), \spa f, \spa F, \spa \mathcal{U}^{(i)}) \Big)^{a_i}.\]
\end{remark}

Suppose $\mathcal{U}^{(i)}\in \mathcal{C}_{X_i}^o$, $i=1,\cdots,r$, we define
\[P^{\boldsymbol{a}}(G,f,\mathcal{U}^{(1)},\cdots,\mathcal{U}^{(r)})=\lim_{n\to \infty}\displaystyle\frac{\log P^{\boldsymbol{a}}(X_r, \spa f,\spa F_n, \spa \mathcal{U}^{(r)})}{|F_n|}.\]
Since $F\in \mathcal{F}(G)\mapsto \log P^{\boldsymbol{a}}(X_r, \spa f,\spa F, \spa \mathcal{U}^{(r)})$ is sub-additive, the limit exists.

Let $\mu\in \mathcal{M}^G(X_1)$. We put
\[P^{\boldsymbol{a}}_{var}(G,f,\mathcal{U}^{(1)},\cdots,\mathcal{U}^{(r)})=\sup_{\mu\in\mathcal{M}^G(X_1)}\left(\sum_{i=1}^rw_ih_{\mu_i}(\mathcal{U}^{(i)},G)+\omega_1\int_X fd\mu\right).\]
\begin{claim}\label{cl3.5}
    Let $\mathcal{U}^{(i)}\in \mathcal{C}_{X_i}^o$. Then for any $\mu\in \mathcal{M}^G(X_1)$,
    \[\sum_{i=1}^rw_ih_{\mu_i}(\mathcal{U}^{(i)},G)+\omega_1\int_X fd\mu\leq P^{\boldsymbol{a}}(G,f,\mathcal{U}^{(1)},\cdots,\mathcal{U}^{(r)}),\]
    that is,
    \[P^{\boldsymbol{a}}_{var}(G,f,\mathcal{U}^{(1)},\cdots,\mathcal{U}^{(r)})\leq P^{\boldsymbol{a}}(G,f,\mathcal{U}^{(1)},\cdots,\mathcal{U}^{(r)}).\]
\end{claim}
\begin{proof}
    Let $\{F_n\}_{n\in \mathbb{N}}$ be a F$\o$lner sequence in $G$ and $\vep>0$. For each $n\in \mathbb{N}$, by Remark \ref{remark3.4}, we can find a partition $\{V_{r1},\cdots,V_{rn_r}\}=\mathcal{V}_r\in \mathcal{P}_{X_r}$ finer than $\mathcal{U}^{(r)}_{F_n}$ such that
    \[\sum_{j_r=1}^{n_r} \Big( P^{\boldsymbol{a}}_{r-1}(\pi_{r-1}^{-1}(V_{rj_r}), \spa f, \spa F_n, \spa \mathcal{U}^{(r-1)}) \Big)^{a_{r-1}}\leq P^{\boldsymbol{a}}(X_r, \spa f,\spa F_n, \spa \mathcal{U}^{(r)})+\vep/2.\]
    Again, by Remark \ref{remark3.4}, we can find a partition $\{V_{r-11},\cdots,V_{r-1n_{r-1}}\}=\mathcal{V}_{r-1}\in \mathcal{P}_{X_{r-1}}$ finer than $\pi_{r-1}^{-1}(\mathcal{V}_r)\cap \mathcal{U}^{(r-1)}_{F_n}$ such that
    \begin{align*}
        &\sum_{j_{r}=1}^{n_{r}}\Big(\sum_{j_{r-1}=1}^{n_{r-1}}\Big(P^{\boldsymbol{a}}_{r-2}(\pi_{r-2}^{-1}(\pi_{r-1}^{-1}(V_{rj_r})\cap V_{r-1j_{r-1}}), \spa f, \spa F_n, \spa \mathcal{U}^{(r-2)})\Big)^{a_{r-2}}\Big)^{a_{r-1}}\\&\leq\sum_{j_r=1}^{n_r} \Big( P^{\boldsymbol{a}}_r(\pi_r^{-1}(V_{rj_r}), \spa f, \spa F_n, \spa \mathcal{U}^{(r-1)}) \Big)^{a_{r-1}}+\vep/2^2.
    \end{align*}
    Hence,
    \begin{align}\label{eq3.3}
        &\sum_{j_{r}=1}^{n_{r}}\Big(\sum_{j_{r-1}=1}^{n_{r-1}}\Big(P^{\boldsymbol{a}}_{(r-2)}(\pi_{r-2}^{-1}(\pi_{r-1}^{-1}(V_{rj_r})\cap V_{r-1j_{r-1}}), \spa f, \spa F_n, \spa \mathcal{U}^{(r-2)})\Big)^{a_{r-2}}\Big)^{a_{r-1}}\\&\leq P^{\boldsymbol{a}}(X_r, \spa f,\spa F_n, \spa \mathcal{U}^{(r)})+\vep/2+\vep/2^2.\notag
    \end{align}
    If $\pi_{r-1}^{-1}(V_{rj_r})\cap V_{r-1j_{r-1}}=\emptyset$, we let
    \[P^{\boldsymbol{a}}_{r-2}\big(\pi_{r-2}^{-1}(\pi_{r-1}^{-1}(V_{rj_r})\cap V_{r-1j_{r-1}}), \spa f, \spa F_n, \spa \mathcal{U}^{(r-2)}\big)^{a_{r-2}}=0.\]
    Since $\mathcal{V}_{r-1}$ is finer than $\pi^{-1}_{r-1}(\mathcal{V}_r)$, \eqref{eq3.3} has the form
    \[\sum_{V_r\in \mathcal{V}_r}\Big(\sum_{\pi_{r-1}V_{r-1}\subset V_r}\Big(P^{\boldsymbol{a}}_{r-2}(\pi^{-1}_{r-2}(V_{r-1}), \spa f, \spa F_n, \spa \mathcal{U}^{(r-2)})\Big)^{a_{r-2}}\Big)^{a_{r-1}}.\]

    By induction, we can find a finite collection of partitions $\{\mathcal{V}_{i}\in \mathcal{P}_{X_i}:i=1,\cdots,r\}$ such that
    \begin{enumerate}
        \item $\mathcal{V}_{i}$ is finer than $\mathcal{U}^{(i)}_{F_n}\cap\pi^{-1}_{i}\mathcal{V}_{i+1}$.
        \item For any $1\leq i< r$, we have
        \begin{align*}
            &\sum_{V_r\in \mathcal{V}_r}\Big(\sum_{\pi_{r-1}V_{r-1}\subset V_r}\Big(\cdots\sum_{\pi_{i+1}V_{i+1}\subset V_{i+2}}\Big(P_{i}^{\boldsymbol{a}}(\pi^{-1}_{i}(V_{i+1}),f,F_n,\mathcal{U}^{(i)})\Big)^{a_{i}}\cdots\Big)^{a_{r-2}}\Big)^{a_{r-1}}\\&\leq\sum_{V_r\in \mathcal{V}_r}\Big(\sum_{\pi_{r-1}V_{r-1}\subset V_r}\Big(\cdots\sum_{\pi_{i+2}V_{i+2}\subset V_{i+3}}\Big(P_{i+1}^{\boldsymbol{a}}(\pi^{-1}_{i+1}(V_{i+2}),f,F_n,\mathcal{U}^{(i+1)})\Big)^{a_{i+1}}\cdots\Big)^{a_{r-2}}\Big)^{a_{r-1}}+\vep/2^{r-i}.
        \end{align*}
    \end{enumerate}
      Hence, for all $i=1,\cdots,r$, we have the inequality
    \begin{align*}
        &\sum_{V_r\in \mathcal{V}_r}\Big(\sum_{\pi_{r-1}V_{r-1}\subset V_{r}}\Big(\cdots\sum_{\pi_{i+1}V_{i+1}\subset V_{i+2}}\Big(P_{i}^{\boldsymbol{a}}(\pi^{-1}_{i}(V_{i+1}),f,F_n,\mathcal{U}^{(i)})\Big)^{a_{i}}\cdots\Big)^{a_{r-2}}\Big)^{a_{r-1}}\\&\leq P^{\boldsymbol{a}}(X_r, \spa f,\spa F_n, \spa \mathcal{U}^{(r)})+\vep.
    \end{align*}
    \begin{remark}\label{rmk3.6}
       It is easy to see that by induction we obtain
        \begin{align}\label{e3.4}
           &\sum_{V_r\in \mathcal{V}_r}\Big(\sum_{\pi_{r-1}V_{r-1}\subset {V}_{r}}\Big(\cdots\sum_{V_{1}\subset\pi_1^{-1}(V_2)}\Big(e^{\sup_{x\in V_1}S_{F_n}f(x)}\Big)^{a_{1}}\cdots\Big)^{a_{r-2}}\Big)^{a_{r-1}}\notag\\&\leq P^{\boldsymbol{a}}(X_r, \spa f,\spa F_n, \spa \mathcal{U}^{(r)})+\vep.
        \end{align}
    \end{remark}

  Let $\mu\in \mathcal{M}^G(X_1)$. Reviewing the definition of $h_\mu(\mathcal{U})$ and by Lemma \ref{lem2.7}, we obtain
    \begin{align*}
&\sum_{i=1}^rw_ih_{\mu_i}(\mathcal{U}^{(i)})+w_1\int_{X_1}fd{\mu}\\&=\lim_{n\to \infty}\displaystyle\frac{1}{|F_n|}\Big(\sum_{i=1}^rw_iH_{\mu_i}(\mathcal{U}^{(i)}_{F_n})+a_1\cdots a_{r-1}|F_n|\int_{X_1}fd{\mu}\Big)\\&\leq\lim_{n\to \infty}\displaystyle\frac{1}{|F_n|}\Big(\sum_{i=1}^rw_iH_{\mu_i}(\mathcal{V}_i)+a_1\cdots a_{r-1}|F_n|\int_{X_1}fd{\mu}\Big),
    \end{align*}
where $\{\mathcal{V}_i\succeq \mathcal{U}^{(i)}_{F_n}:i=1,\cdots,r\}$ are partitions generated via the above procedure satisfying conditions (1) and (2). From the definition of $\boldsymbol{\omega}$, we have
\begin{align}\label{e3.5}
    &\sum_{i=1}^rw_iH_{\mu_i}(\mathcal{V}_i)+a_1\cdots a_{r-1}|F_n|\int_{X_1}fd{\mu}\notag\\&=H_{\mu_r}(\mathcal{V}_r)+a_1\cdots a_{r-1}\int_{X_1}S_{F_n}fd{\mu}+\sum_{i=1}^{r-1}a_ia_{i+1}\cdots a_{r-1}\Big(H_{\mu_i}(\mathcal{V}_i)-H_{\mu_{i+1}}(\mathcal{V}_{i+1})\Big)\notag\\&=
    H_{\mu_r}(\mathcal{V}_r)+a_1\cdots a_{r-1}\int_{X_1}S_{F_n}fd{\mu}+\sum_{i=1}^{r-1}a_ia_{i+1}\cdots a_{r-1}\Big(H_{\mu_i}(\mathcal{V}_i|\pi_{i}^{-1}\mathcal{V}_{i+1})\Big).
\end{align}
Here, we use the fact that
$\pi^{-1}_i\mathcal{V}_{i+1}\preceq\mathcal{V}_i \text{ (due to (1))}.$
Fix $F_n$, we use backward induction to evaluate \eqref{e3.5}. First, consider the term
\begin{align*}
    &a_1\cdots a_{r-1}\int_{X_1}S_{F_n}fd{\mu}+a_1a\cdots a_{r-1}H_{\mu}(\mathcal{V}_1|\pi_{1}^{-1}\mathcal{V}_{2})\\&=a_1\cdots a_{r-1}\Big(\int_{X_1}S_{F_n}fd{\mu}+H_{\mu}(\mathcal{V}_1|\pi_{1}^{-1}\mathcal{V}_{2})\Big).
\end{align*}
By Lemma \ref{lem3.3}, we have the following inequality:
\begin{align*}
    &\int_{X_1}S_{F_n}fd{\mu}+H_{\mu}(\mathcal{V}_1|\pi_{1}^{-1}\mathcal{V}_{2})\\&\leq\sum_{C\in \mathcal{V}_2}\mu(\pi_1^{-1}C)\Big\{\sum_{D\in \mathcal{V}_1}\Big(-\displaystyle\frac{\mu(D\cap \pi_1^{-1}C)}{\mu(\pi_1^{-1}C)}\log\displaystyle\frac{\mu(D\cap \pi_1^{-1}C)}{\mu(\pi_1^{-1}C)}+\displaystyle\frac{\mu(D\cap \pi_1^{-1}C)}{\mu(\pi_1^{-1}C)}\sup_{x\in D}S_{F_n}f(x)\Big)\Big\}\\&\leq\sum_{C\in \mathcal{V}_2}\mu(\pi_1^{-1}C)\log\sum_{D\subset \pi_1^{-1}C}e^{\sup_{x\in D}S_{F_n}f(x)}.
\end{align*}
Applying this inequality to \eqref{e3.4}, we have
\begin{align*}
    &a_2a_3\cdots a_{r-1}\Big(a_1\int_{X_1}S_{F_n}fd{\mu}+a_1H_{\mu}(\mathcal{V}_1|\pi_{1}^{-1}\mathcal{V}_{2})+H_{\mu_2}(\mathcal{V}_2|\pi_{2}^{-1}\mathcal{V}_{3})\Big)\notag\\&\leq a_2a_3\cdots a_{r-1}\Big(a_1\sum_{C\in \mathcal{V}_2}\mu(\pi_1^{-1}C)\log\sum_{D\subset \pi_1^{-1}C}e^{\sup_{x\in D}S_{F_n}f(x)}+H_{\mu_2}(\mathcal{V}_2|\pi_{2}^{-1}\mathcal{V}_{3})\Big).
\end{align*}
Using Lemma \ref{lem3.3} again, we have
\begin{align*}
    &a_1\sum_{C\in \mathcal{V}_2}\mu(\pi_1^{-1}C)\log\sum_{D\subset \pi_1^{-1}C}e^{\sup_{x\in D}S_{F_n}f(x)}+H_{\mu_2}(\mathcal{V}_2|\pi_{2}^{-1}\mathcal{V}_{3})\\&=\sum_{C\in \mathcal{V}_3}\mu_3(C)\Big\{\sum_{D\in \mathcal{V}_2}\Big(-\displaystyle\frac{\mu_2(D)}{\mu_3(C)}\log \displaystyle\frac{\mu_2(D)}{\mu_3(C)}+\displaystyle\frac{\mu_2(D)}{\mu_3(C)}\log\Big(\sum_{E\subset \pi_1^{-1}D}e^{\sup_{x\in D}S_{f_n(x)}}\Big)^{a_1}\Big)\Big\}\\&\leq\sum_{C\in \mathcal{V}_3}\mu_3(C)\log\sum_{D\in \mathcal{V}_2}\Big(\sum_{E\subset \pi_1^{-1}D}e^{\sup_{x\in D}S_{f_n(x)}}\Big)^{a_1}.
\end{align*}
Continue likewise, we obtain the following upper bound of \eqref{e3.5}:
\begin{align}
     \log\sum_{C_r\in \mathcal{V}_r}\Big(\sum_{\pi_{r-1}C_{r-1}\subset C_r}\Big(\cdots\sum_{C_{1}\subset\pi^{-1}(C_2)}\Big(e^{\sup_{x\in C_1}S_{F_n}f(x)}\Big)^{a_{1}}\cdots\Big)^{a_{r-2}}\Big)^{a_{r-1}}.
\end{align}
Hence, from Remark \ref{rmk3.6}, we have
\begin{align*}
    &\sum_{i=1}^rw_iH_{\mu_i}(\mathcal{V}_i)+a_1\cdots a_{r-1}|F_n|\int_{X_1}fd{\mu}\\&\leq \log\sum_{C_r\in \mathcal{V}_r}\Big(\sum_{C_{r-1}\in \mathcal{V}_{r-1}}\Big(\cdots\sum_{C_{1}\subset\pi_1^{-1}(C_2)}\Big(e^{\sup_{x\in C_1}S_{F_n}f(x)}\Big)^{a_{1}}\cdots\Big)^{a_{r-2}}\Big)^{a_{r-1}}\\&\leq \log \big(P^{\boldsymbol{a}}(X_r, \spa f,\spa F_n, \spa \mathcal{U}^{(r)})+\vep\big).
\end{align*}
Since $\mathcal{V}_i$ is finer than $\mathcal{U}^{(i)}_{F_n}$, we obtain
\begin{align*}
    &\displaystyle\frac{1}{|F_n|}\sum_{i=1}^rw_iH_{\mu_i}(\mathcal{U}^{(i)}_{F_n})+a_1\cdots a_{r-1}|F_n|\int_{X_1}fd{\mu}\\&\leq\displaystyle\frac{1}{|F_n|}\sum_{i=1}^rw_iH_{\mu_i}(\mathcal{V}_i)+a_1\cdots a_{r-1}|F_n|\int_{X_1}fd{\mu}\\&\leq\displaystyle\frac{1}{|F_n|}\log \big(P^{\boldsymbol{a}}(X_r, \spa f,\spa F_n, \spa \mathcal{U}^{(r)})+\vep\big).
\end{align*}
Since $\vep$ can be chosen arbitrarily, and by taking the limits on both sides, we get
\begin{align*}
    P^{\boldsymbol{a}}_{var}(G,f,\mathcal{U}^{(1)},\cdots,\mathcal{U}^{(r)})\leq P^{\boldsymbol{a}}(G,f,\mathcal{U}^{(1)},\cdots,\mathcal{U}^{(r)}).
\end{align*}
\end{proof}
\begin{claim}\label{cl3.7}
    For any $\mathcal{U}^{(i)}\in \mathcal{C}^o_{X_i}$, $i=1,\cdots,i$, we have
    \[P^{\boldsymbol{a}}(G,f,\mathcal{U}^{(1)},\cdots,\mathcal{U}^{(r)})\leq P^{\textbf{a}}(f,G).\]
\end{claim}
\begin{proof}
Let $\delta_i$ be the Lebesgue number of $\mathcal{U}^{(i)}$, that is, if $K_i$ is a subset of $X$ such that its diameter is less than $\delta_i$, then $K_i\subset U_i$ for some $U_i\in \mathcal{U}^{(i)}$. Let $\vep$ be a positive number with $\vep<\underset{i}{\min}\,\{\delta_i:i=1,\cdots,r\}$. Notice that if $\mathcal{V}_i$ is an open cover such that $\diam(\mathcal{V}_i,d_F^{(i)})<\vep$ with respect to $F\in \mathcal{F}(G)$, then for any $x,y\in V\in \mathcal{V}_i$, we have $d^{(i)}(gx,gy)<\vep$. Hence, there is a $U\in \bigvee_{g\in F}g^{-1}\mathcal{U}^{(i)}$ such that $x,y\in U$, which means that $V\subset U$ and $\mathcal{V}_i\succeq\mathcal{U}^{(i)}_F$. Therefore,
\[{P}^{\textbf{a}}(X_r,f,F,\epsilon)\geq P^{\boldsymbol{a}}(X_r, \spa f,\spa F, \spa \mathcal{U}^{(r)}),\]
for all $F\in \mathcal{F}(G)$. Hence,
\[P^{\textbf{a}}(f,G)\geq P^{\boldsymbol{a}}(G,f,\mathcal{U}^{(1)},\cdots,\mathcal{U}^{(r)}).\]
\end{proof}
 Theorem \ref{thm3.1} follows from Remark \ref{rmk2.9} (3) and Claims \ref{cl3.5} and \ref{cl3.7}.
\end{proof}
\section{Proof of $P_{var}^{\textbf{a}}(G,f)\geq P^{\textbf{a}}(G,f)$ and $h_{var}^{\textbf{a}}(G)\geq h_{top}^{\textbf{a}}(G)$}
In this section, we aim to prove the inverse direction of the variational principle as followed.
\begin{theorem}\label{thm4.1}
    Let $(X_i, G)$, $i=1,\cdots,r$, be zero-dimensional dynamical systems with discrete countable amenable group action and $\pi_{i}: X_i\to X_{i+1}$, $i=1,\cdots,r-1$ be factor maps. Then we obtain the \textit{lower bound}
    \[P_{var}^{\textbf{a}}(G,f,\boldsymbol{\pi})\geq P^{\textbf{a}}(G,f,\boldsymbol{\pi}),\]
    for any $f\in C(X)$.
\end{theorem}
The following result is a direct consequence of Theorem \ref{thm3.1}.
\begin{corollary}
     Let $(X_i, G)$, $i=1,\cdots,r$, be zero-dimensional dynamical systems with discrete countable amenable group action and $\pi_{i}: X_i\to X_{i+1}$, $i=1,\cdots,r-1$ be factor maps. Then we obtain the \textit{lower bound}
     \[h_{var}^{\textbf{a}}(G,\boldsymbol{\pi})\geq h^{\textbf{a}}(G,\boldsymbol{\pi}).\]
\end{corollary}

\begin{proof}[Proof of Theorem \ref{thm4.1}]
Let $\mathcal{U}^{(i)}\in \mathcal{C}^o_{X_i}$ be clopen covers of $X_i$ such that $\pi^{-1}_{i-1}(\mathcal{U}^{(i)})\preceq\mathcal{U}^{(i-1)}$. We define
\begin{flalign*}
& Q^{\boldsymbol{a}}_1(X_1, \spa f, \spa F, \spa \mathcal{U}^{(1)}) &
\end{flalign*} \\[-35pt]
\begin{align}
&= \inf \left\{ \spa \sum_{j=1}^n \exp \left( \sup_{U_j} S_F f \right) \spa \middle|
\begin{array}{l}
\text{$n \in \mathbb{N}$, $\{V_j\}_{j=1}^n$ is a cover of $X_1$} \\
\text{with $V_i\in\mathcal{V}\succeq\mathcal{U}^{(1)}_F$  for some $\mathcal{V}\in \mathcal{C}_{X_1}$ }
\end{array}
\right\}.
\end{align}\label{e4.1}
Suppose $Q_i^{\textbf{a}}(X_i, \spa f, \spa F, \spa \mathcal{U}^{(i)})$ is already defined. We define
\begin{flalign*}
& Q^{\boldsymbol{a}}_{i+1}(X_{i+1}, \spa f,\spa F, \spa \mathcal{U}^{(i+1)}) &
\end{flalign*} \\[-35pt]
\begin{align}\label{e4.2}
&= \inf \left\{ \spa \sum_{j=1}^n \Big( Q^{\boldsymbol{a}}_i(\pi_i^{-1}(V_j), \spa f, \spa F, \spa \mathcal{U}^{(i)}) \Big)^{a_i} \spa \middle|
\begin{array}{l}
\text{$n \in \mathbb{N}$, $\{V_j\}_{j=1}^n$ is a cover of $X_{i+1}$} \\
\text{with $V_j\in \mathcal{V}\succeq\mathcal{U}^{(i+1)}_F$ for some $\mathcal{V}\in \mathcal{C}_{X_{i+1}}$}
\end{array}
\right\}.
\end{align}
\begin{remark}
   Similar to Remark \ref{remark3.4}, the cover $\mathcal{V}$ in \eqref{e4.2} can be chosen from partitions. Since $\mathcal{U}^{(i)}$ $(i=1,\cdots,r)$ are clopen covers and $\pi^{-1}_{i-1}(\mathcal{U}^{(i)})\preceq\mathcal{U}^{(i-1)}$, it is easy to check that the cover $\mathcal{V}$ in \eqref{e4.2} is identitical with $\mathcal{U}^{(i+1)}$.
\end{remark}

We define
\[Q^{\boldsymbol{a}}(G,f,\mathcal{U}^{(1)},\cdots,\mathcal{U}^{(r)})=\lim_{n\to \infty}\displaystyle\frac{\log Q^{\boldsymbol{a}}(X_r, \spa f,\spa F_n, \spa \mathcal{U}^{(r)})}{|F_n|}.\]
Since $F\in \mathcal{F}(G)\mapsto \log Q^{\boldsymbol{a}}(X_r, \spa f,\spa F, \spa \mathcal{U}^{(r)})$ is sub-additive, the limit exists. Then we define
\[Q^{\boldsymbol{a}}(G,f,\boldsymbol{\pi})=\lim_{diam(\mathcal{U}^{(i)})\to 0}Q^{\boldsymbol{a}}(G,f,\mathcal{U}^{(1)},\cdots,\mathcal{U}^{(r)}),\]
where $i=1,\cdots,r$ and $\mathcal{U}^{(i)}\in \mathcal{C}_X^o$ are clopen partitions of $X_i$ with $\pi_{i-1}^{-1}(\mathcal{U}^{(i)})\preceq\mathcal{U}^{(i-1)}$.
\begin{claim}\label{cl4.3}
   Under the above conditions, we have
    \[P^{\boldsymbol{a}}(G,f,\boldsymbol{\pi})=Q^{\boldsymbol{a}}(G,f,\boldsymbol{\pi}).\]
\end{claim}
\begin{proof}
    From the definition of ${P}^{\textbf{a}}(X_r,f,\vep)$, if $diam(\mathcal{U}_i)\leq \vep$, then for every $F\in \mathcal{F}(G)$, $\diam(\mathcal{U}^{(i)}_F)\leq \vep$ and
    \[Q^{\boldsymbol{a}}(G,f,\mathcal{U}^{(1)},\cdots,\mathcal{U}^{(r)})\geq {P}^{\textbf{a}}(X_r,f,\vep).\]
    Therefore,
    \[Q^{\boldsymbol{a}}(G,f,\boldsymbol{\pi})\geq P^{\textbf{a}}(G,f,\boldsymbol{\pi}).\]
For any clopen partitions $\mathcal{U}^{(i)}\in \mathcal{C}_X^o$ with $\pi_{i-1}^{-1}(\mathcal{U}^{(i)})\preceq\mathcal{U}^{(i)}$, let $\delta$ be the minimal Lebesgue number of $\mathcal{U}^{(i)}$ $(i=1,\cdots,r)$. Then we have
\[Q^{\boldsymbol{a}}(G,f,\mathcal{U}^{(1)},\cdots,\mathcal{U}^{(r)})\leq {P}^{\textbf{a}}(X_r,f,\delta),\]
hence,
\[P^{\textbf{a}}(G,f,\boldsymbol{\pi})=Q^{\boldsymbol{a}}(G,f,\boldsymbol{\pi}).\]
\end{proof}
Now we shall prove the following statement.
\begin{lemma}
   For any clopen partitions $\mathcal{U}^{(i)}\in \mathcal{C}_X^o$ with $\pi_{i-1}^{-1}(\mathcal{U}^{(i)})\preceq\mathcal{U}^{(i-1)}$ $i=1,\cdots,r$, there is a $\mu\in \mathcal{M}^G(X_1)$ such that
    \begin{align}
        Q^{\boldsymbol{a}}(G,f,\mathcal{U}^{(1)},\cdots,\mathcal{U}^{(r)})\leq \sum_{i=1}^rw_ih_{\mu_i}(\mathcal{U}^{(i)},G)+\omega_1\int_X fd\mu,
    \end{align}
     where $\mu_i=\pi^{(i-1)}\mu$.
\end{lemma}
\begin{proof}
    Let $\{F_n:n\in \mathbb{N}\}$ be a F$\o$lner sequence in $G$. For any $i<j$ and $U\in \mathcal{U}^{(i)}_{F_n}$, we put
    \[\mathcal{U}^{(i)}_{F_n}(U)=\{V\in \mathcal{U}^{(i)}_{F_n}:\pi_{j-1}\circ\pi_{j-2}\circ\cdots\circ\pi_{i}(V)\subset U\}.\]
If $V\in \mathcal{U}^{(i)}_{F_n}$ and $i<j$, we denote
\[\tpi_j(V)=U,\]
the unique element in $\mathcal{U}^{(j)}_{F_n}$ such that $\pi_{j-1}\circ\cdots\circ\pi_{i}(V)\subset U$.

    For any $U\in \mathcal{U}^{(2)}_{F_n}$, we define a map $Z^{(1)}_{F_n}$ from $\mathcal{U}^{(2)}_{F_n}$ to $\mathbb{R}$ by
\begin{align*}
    Z^{(1)}_{F_n}(U)=\sum_{V\in \mathcal{U}^{(1)}_{F_n}(U)}e^{\sup_{x\in V}S_{F_n}f(x)}.
\end{align*}
Suppose the map $Z^{(i-1)}_{F_n}$ from $\mathcal{U}^{(i)}_{F_n}$ to $\mathbb{R}$ has been defined. For any $U\in \mathcal{U}^{(i+1)}_{F_n}$, we set
\begin{align*}
    Z^{(i)}_{F_n}(U)=\sum_{V\in \mathcal{U}^{(i)}_{F_n}(U)}\Big(Z^{(i-1)}_{F_n}(V)\Big)^{a_{i-1}}.
\end{align*}
Then we define
\begin{align*}
    Z_{F_n}=\sum_{U\in \mathcal{U}^{(r)}_{F_n}}\Big(Z^{(r-1)}_{F_n}(V)\Big)^{a_{r-1}}.
\end{align*}
Actually, from the definition of $Q^{\boldsymbol{a}}$, we have
\begin{align*}
    Z_{F_n}=Q^{\boldsymbol{a}}(X_r, \spa f,\spa F_n, \spa \mathcal{U}^{(r)}).
\end{align*}

 Hence, it suffices to show that there is an invariant measure $\mu\in \mathcal{M}^G(X_1)$ such that
\[\sum_{i=1}^rw_ih_{\mu_i}(\mathcal{U}_i,G)+\omega_1\int_X fd\mu\leq\lim_{n\to \infty}\displaystyle\frac{\log Z_{F_n}}{|F_n|}.\]

We notice that it is a similar way in \cite{al} to construct such invariant measure for discrete amenable group $G$.

Since the elements in $\mathcal{U}^{(1)}_{F_n}$ are clopen, for each $V\in \mathcal{U}^{(1)}_{F_n}$ we can find a point $x_V$ such that
\[e^{S_{F_n}f(x_V)}=\sup_{x\in V}e^{S_{F_n}f(x)}.\]
Fix $F_n$ we define a probability measure on $X_1$ by
\begin{align*}
\nu_{n} =
\begin{multlined}[t][11.5cm]
\frac{1}{Z_{F_n}} \sum_{V \in \mathcal{U}^{(1)}_{F_n}} {Z_{F_n}^{(r-1)}(\tpi_{r-1}(V))}^{a_{r-1}-1}{Z_{F_n}^{(r-2)}(\tpi_{r-2}(V))}^{a_{r-2}-1} \\
\times \cdots \times {Z_{F_n}^{(2)}(\tpi_3 (V))}^{a_{2}-1}{Z_{F_n}^{(1)}(\tpi_2 (V))}^{a_{1}-1}e^{S_{F_n} f(x_V)} \delta_{x_V},
\end{multlined}
\end{align*}
where $\delta_{x_V}$ is the Dirac measure at $x_A$. Moreover, $\nu_{n}$ is a probability measure on $X_1$, that is, $\nu_n(X_1)=1$.
In general, $\nu_{n}$ is not a $G$-invariant measure. Hence, we set
\[\mu_n=\displaystyle\frac{1}{|F_n|}\sum_{g\in F_n}g\nu_{n},\]
where $g\nu_n(\cdot)=\nu_n(g^{-1}\cdot)$. By passing to subsequence twice in succession, $\{\mu_n\}$ convergences to an invariant probability measure $\mu$ on $X_1$ in the compact weak*-topology of $\mathcal{M}(X_1)$.

We arrange $\{F_n\}_{n\in \mathbb{N}}$ such that
\begin{enumerate}
    \item the limit $\lim_{n\to \infty}\log Z_{F_n}/{|F_n|}$ exists, and
    \item the sequence $\{\mu_n\}$ convergence to some $\mu\in \mathcal{M}^G(X_1)$.
\end{enumerate}
For convenience, we denote
\begin{align*}
    \nu_n^{(i)}&=\pi^{(i-1)}\nu_n\\&=\frac{1}{Z_{F_n}} \sum_{V \in \mathcal{U}^{(1)}_{F_n}} {Z_{F_n}^{(r-1)}(\tpi_{r}(V))}^{a_{r-1}-1} \cdots {Z_{F_n}^{(1)}(\tpi_2 (V))}^{a_{1}-1}e^{S_{F_n} f(x_V)} \delta_{\pi^{(i-1)}(x_V)},
\end{align*}
where $2\leq i\leq r$, and write
\begin{align}\label{e4}
    W_{F_n}^{(j)}=\sum_{U\in \mathcal{U}_{F_n}^{(j+1)}}{Z_{F_n}^{(r-1)}(\tpi_r U)}^{a_{r-1}-1}
\cdots {Z_{F_n}^{(j+1)}(\tpi_{j+2} U)}^{a_{j+1}-1} {Z_{F_n}^{(j)}(U)}^{a_{j}} \log{\left(Z_{F_n}^{(j)} (U) \right)}.
\end{align}
\begin{claim}\label{cl4.6}
   Suppose $\nu_n$ are probability measure defined above, we obtain the following equalities:
   \[H_{\nu_n}(\mathcal{U}_{F_n}^{(1)})=\log Z_{F_n}-\int_{X_1}S_{F_n}fd\nu_n-\sum_{j-1}^{r-1}\frac{a_j-1}{Z_{F_n}}W_{F_n}^{(j)},\]
   and
   \[H_{\nu_n^{(i)}}(\mathcal{U}^{(i)}_{F_n})=\log Z_{F_n}-\frac{a_{i-1}}{Z_{F_n}}W_{F_n}^{(i-1)}-\sum_{j=i}^{r-1}\frac{a_j-1}{Z_{F_n}}W_{F_n}^{(j)},\]
   where $i=2,\cdots,r$.
\end{claim}
\begin{proof}
    For any $V\in \mathcal{U}_{F_n}^{(1)}$, we have
    \begin{align*}
        \nu_n(V)=\frac{1}{Z_{F_n}} {Z_{F_n}^{(r-1)}(\tpi_{r}(V))}^{a_{r-1}-1} \cdots {Z_{F_n}^{(1)}(\tpi_2 (V))}^{a_{1}-1}e^{S_{F_n} f(x_V)}.
    \end{align*}
    Hence,
    \begin{flalign*}
& H_{\nu_n}(\mathcal{U}_{F_n}^{(1)}) = - \sum_{V \in\mathcal{U}_{F_n}^{(1)}} \nu_n(V) \log{(\nu_n(V))} &
\end{flalign*} \\[-30pt]
\begin{flalign}\label{e4.4}
&
\begin{multlined}[t][11.5cm]
= \spa \log{Z_{F_n}} - \underbrace{\sum_{V \in \mathcal{U}_{F_n}^{(1)}} \nu_n(V) S_{F_n} f(x_V)}_{(\mathrm{I})} \\
- \sum_{j=1}^{r-1} \frac{a_j-1}{Z_{F_n}} \underbrace{\sum_{V \in \mathcal{U}_{F_n}^{(1)}} {Z_{F_n}^{(r-1)}(\tpi_r V)}^{a_{r-1}-1} \cdots {Z_{F_n}^{(1)}(\tpi_2 V)}^{a_{1}-1} e^{S_{F_n} f(x_V)} \log{\left(Z_{F_n}^{(j)}(\tpi_{j+1} V) \right)}}_{(\mathrm{I}\hspace{-0.5pt}\mathrm{I})}.
\end{multlined}
\end{flalign}
We calculate $\mathrm{(I)}$ by
\begin{align}\label{e4.5}
    \mathrm{(I)}&=\int_{X_1}S_{F_n}fd\nu_n=\sum_{V \in \mathcal{U}_{F_n}^{(1)}} \nu_n(V) S_{F_n} f(x_V)\notag\\&=\sum_{V\in \mathcal{U}_{F_n}^{(1)}}\frac{1}{Z_{F_n}} {Z_{F_n}^{(r-1)}(\tpi_{r}(V))}^{a_{r-1}-1} \cdots {Z_{F_n}^{(1)}(\tpi_2 (V))}^{a_{1}-1}e^{S_{F_n} f(x_V)}S_{F_n} f(x_V).
\end{align}
Recall that for any $V\in \mathcal{U}_{F_n}^{(1)}$ and $j>0$, there is a unique clopen set $U\in \mathcal{U}_{F_n}^{(j+1)}$ such that $\tpi_{j+1}(V)=U$. Hence,
\begin{align*}
    \mathrm{(II)}&=\sum_{U\in \mathcal{U}_{F_n}^{(j+1)}}\sum_{V\in \mathcal{U}_{F_n}^{(1)}(U)} {Z_{F_n}^{(r-1)}(\tpi_r V)}^{a_{r-1}-1} \cdots {Z_{F_n}^{(1)}(\tpi_2 V)}^{a_{1}-1} e^{S_{F_n} f(x_V)} \log{\left(Z_{F_n}^{(j)}(\tpi_{j+1} V) \right)}\\
&=
\begin{multlined}[t][8cm]
\sum_{U \in \mathcal{U}_{F_n}^{(j+1)}} {Z_{F_n}^{(r-1)}(\tpi_r U)}^{a_{r-1}-1}
\cdots {Z_{F_n}^{(j+1)}(\tpi_{j+2} U)}^{a_{j+1}-1} {Z_{F_n}^{(j)}(U)}^{a_{j}-1} \log{\left(Z_{F_n}^{(j)} (U) \right)} \\
\times \left(\underbrace{\sum_{V \in \mathcal{U}_{F_n}^{(1)}(U)} {Z_{F_n}^{(j-1)}(\tpi_j V)}^{a_{j-1}-1} \cdots {Z_{F_n}^{(1)}(\tpi_2 V)}^{a_{1}-1}e^{S_{F_n} f(x_V)}}_{(\mathrm{I}\hspace{-0.5pt}\mathrm{I})'}\right),
\end{multlined}
\end{align*}
where $\mathrm{(II)'}$ can be calculated by
\begin{align*}
    &\mathrm{(II)'}=\sum_{V \in \mathcal{U}_{F_n}^{(1)}(U)} {Z_{F_n}^{(j-1)}(\tpi_j V)}^{a_{j-1}-1} \cdots {Z_{F_n}^{(1)}(\tpi_2 V)}^{a_{1}-1}e^{S_{F_n} f(x_V)}\\&=\begin{multlined}[t][15cm]
\sum_{V_j \in \mathcal{U}_{F_n}^{(j)}(U)} {Z_{F_n}^{(j-1)}(V_j)}^{a_{j-1}-1} \hspace{-10pt} \sum_{V_{j-1} \in \mathcal{U}_{F_n}^{(j-1)}(V_{j})} {Z_{F_n}^{(j-2)}(V_{j-1})}^{a_{j-2}-1} \\
\cdots \sum_{V_3 \in \mathcal{U}_{F_n}^{(3)}(V_4)} Z_{F_n}^{(2)}(V_3)^{a_2-1} \sum_{V_2 \in \mathcal{U}^{(2)}_{F_n}(V_3)} {Z_{F_n}^{(1)}(V_2)}^{a_{1}-1} \underbrace{\sum_{V_1 \in \mathcal{U}_{F_n}^{(1)}(V_2)} e^{S_{F_n} f(x_{V_1})}}_{= Z_{F_n}^{(1)}(V_2)}
\end{multlined}
\\&=\cdots = \sum_{V_j \in \mathcal{U}_{F_n}^{(j)}(U)} {Z_{F_n}^{(j-1)}(V_j)}^{a_{j-1}} = Z^{(j)}_{F_n}(U).
\end{align*}
Therefore, we have
\begin{align*}
    \mathrm{(II)}&=\sum_{U \in \mathcal{U}_{F_n}^{(j+1)}} {Z_{F_n}^{(r-1)}(\tpi_r U)}^{a_{r-1}-1}
\cdots {Z_{F_n}^{(j+1)}(\tpi_{j+2} U)}^{a_{j+1}-1} {Z_{F_n}^{(j)}(U)}^{a_{j}} \log{\left(Z_{F_n}^{(j)} (U) \right)}\\&=W^{(j)}_{F_n}.
\end{align*}
Combining \eqref{e4.4} and \eqref{e4.5}, we get
\[H_{\nu_n}(\mathcal{U}_{F_n}^{(1)})=\log Z_{F_n}-\int_{X_1}S_{F_n}fd\nu_n-\sum_{j-1}^{r-1}\frac{a_j-1}{Z_{F_n}}W_{F_n}^{(j)}.\]

Now, we turn to calculate $H_{\nu_n^{(i)}}(\mathcal{U}_{F_n}^{(i)})$. For any $U\in \mathcal{U}_{F_n}^{(i)}$,
\begin{align*}
    \nu_n^{(i)}(U)=\frac{1}{Z_{F_n}} \sum_{V\in \mathcal{U}_{F_n}^{(1)}(U)} {Z_{F_n}^{(r-1)}(\tpi_{r}(V))}^{a_{r-1}-1} \cdots {Z_{F_n}^{(1)}(\tpi_2 (V))}^{a_{1}-1}e^{S_{F_n} f(x_V)} \delta_{\pi^{(i-1)}(x_V)}.
\end{align*}
Since $\pi^{-1}_{j}\mathcal{U}_{F_n}^{(j+1)}\preceq\mathcal{U}_{F_n}^{(j)}$ for any $j\geq i$, we have
\begin{equation*}
    \tpi_j(V)=\tpi_j(U).
\end{equation*}
Hence,
\begin{align*}
    \nu_n^{(i)}(U)=\begin{multlined}[t][11.5cm]
\frac{1}{Z_{F_n}} {Z_{F_n}^{(r-1)}(\tpi_r U)}^{a_{r-1}-1} \cdots {Z_{F_n}^{(i-1)}(\tpi_i U)}^{a_{i-1}-1} \\
\times \left(\underbrace{\sum_{V \in \mathcal{U}_{F_n}^{(1)}(U)}
{Z_{F_n}^{(i-2)}(\tpi_{i-1} V)}^{a_{i-2}-1} \cdots {Z_{F_n}^{(1)}(\tpi_2 V)}^{a_{1}-1}e^{S_{F_n} f(x_V)}}_{(\mathrm{I}\hspace{-0.5pt}\mathrm{I})''}\right),
\end{multlined}
\end{align*}
where $\mathrm{(II)}''$ can be calculated as $\mathrm{(II)}'$ above, that is,
\begin{align*}
    \mathrm{(II)}''=Z_{F_n}^{(i-1)}(U).
\end{align*}
Therefore,
\[ \nu_n^{(i)}(U)=\frac{1}{Z_{F_n}} {Z_{F_n}^{(r-1)}(\tpi_r U)}^{a_{r-1}-1} \cdots {Z_{F_n}^{(i-1)}(\tpi_i U)}^{a_{i-1}-1}Z_{F_n}^{(i-1)}(U)^{a_{i-1}}.\]
Thus,
\begin{flalign*}
& \hspace{18pt} H_{\nu_n^{(i)}}(\mathcal{U}^{(i)}_{F_n}) = - \sum_{U \in \mathcal{U}^{(i)}_{F_n}} \nu_{n}^{(i)}(U) \log{\nu_n^{(i)}(U)} &
\end{flalign*} \\[-35pt]
\begin{align*}
&=
\begin{multlined}[t][15cm]
\log{Z_{F_n}} - \frac{1}{Z_{F_n}} \sum_{U \in \mathcal{U}_{F_n}^{(i)}} {Z_{F_n}^{(r-1)}(\tpi_r U)}^{a_{r-1}-1} \cdots {Z_{F_n}^{(i)}(\tpi_{i+1} U)}^{a_{i}-1} {Z_{F_n}^{(i-1)}(U)}^{a_{i-1}} \\
\times \log{\left( {Z_{F_n}^{(r-1)}(\tpi_r U)}^{a_{r-1}-1} \cdots {Z_{F_n}^{(i)}(\tpi_{i+1} U)}^{a_{i}-1} {Z_{F_n}^{(i-1)}(U)}^{a_{i-1}} \right)}
\end{multlined} \\[12pt]
&=
\begin{multlined}[t][16cm]
\log{Z_{F_n}} - \frac{a_{i-1}}{Z_{F_n}} \sum_{U \in \mathcal{U}^{(i)}_{F_n}} {Z_{F_n}^{(r-1)}(\tpi_r U)}^{a_{r-1}-1} \cdots {Z_{F_n}^{(i)}(\tpi_{i+1} U)}^{a_{i}-1}
{Z_{F_n}^{(i-1)}(U)}^{a_{i-1}} \log{\left( Z_{F_n}^{(i-1)}(U) \right)} \\
- \sum_{j=i}^{r-1} \frac{a_{j}-1}{Z_{F_n}} \sum_{U \in \mathcal{U}^{(i)}_{F_n}} \underbrace{{Z_{F_n}^{(r-1)}(\tpi_r U)}^{a_{r-1}-1} \cdots {Z_{F_n}^{(i)}(\tpi_{i+1} U)}^{a_{i}-1} {Z_{F_n}^{(i-1)}(U)}^{a_{i-1}} \log{\left( Z_{F_n}^{(j)}(\tpi_{j+1} U) \right)}}_{\mathrm{(III)}}.
\end{multlined}
\end{align*}
For $\mathrm{(III)}$, we have
\begin{align*}
    \mathrm{(III)}&=\begin{multlined}[t][3cm]
\sum_{U_{j+!} \in \mathcal{U}^{(j+1)}_{F_n}} {Z_{F_n}^{(r-1)}(\tpi_r U_{j+1})}^{a_{r-1}-1} \cdots {Z_{F_n}^{(j+1)}(\tpi_{j+2} U_{j+1})}^{a_j-1} {Z_{F_n}^{(j)}(U_{j+1})}^{a_{j-1}-1} \log{\left( Z_{F_n}^{(j)}(U_{j+1}) \right)} \\
\times \sum_{U_j \in \mathcal{U}_{F_n}^{(j)}(U_{j+1})} {Z^{(j-1)}_{F_n}(U_j)}^{a_{j-2}-1} \cdots \sum_{U_{i+1} \in \mathcal{U}^{(i+1)}_{F_n}(U_{i+2})} {Z^{(i)}_{F_n}(U_{i+1})}^{a_{i+1}-1} \underbrace{\sum_{U_i \in \mathcal{U}^{(i)}_{F_n}(U_{i+1})} {Z^{(i-1)}_{F_n}(U_i)}^{a_{i-1}}}_{= Z^{(i)}_{F_n}(U_{i+1})}
\end{multlined} \\
&= \cdots = \sum_{U_{j+1} \in \mathcal{U}^{(j+1)}_{F_n}} {Z_{F_n}^{(r-1)}(\tpi_r U_{j+1})}^{a_{r-1}-1} \cdots {Z_{F_n}^{(j+1)}(\tpi_{j+2} U_{j+1})}^{a_j-1} {Z_{F_n}^{(j)}(U_{j+1})}^{a_{j-1}} \log{\left( Z_{F_n}^{(j)}(U_{j+1}) \right)}\\&=W_{F_n}^{(j)}\text{ (see }\eqref{e4}).
\end{align*}
Hence, we obtain
\begin{align*}
    H_{\nu_n^{(i)}}(\mathcal{U}^{(i)}_{F_n})=\log Z_{F_n}-\frac{a_{i-1}}{Z_{F_n}}W_{F_n}^{(i-1)}-\sum_{j=i}^{r-1}\frac{a_j-1}{Z_{F_n}}W_{F_n}^{(j)}.
\end{align*}
\end{proof}
By Claim \ref{cl4.6}, we have
\begin{align*}
	\sum_{i=1}^r w_i H_{\nu_{n}^{(i)}} (\mathcal{U}^{(i)}_{F_n}) + w_1 \int_{X_1} S_{F_n} f d\nu_n
	&= \spa \log{Z_{F_n}} - \sum_{i=2}^r \frac{w_i a_{i-1}}{Z_{F_n}} W_{F_n}^{(i-1)} - \sum_{i=1}^{r-1} \sum_{j=i}^{r-1} \frac{w_i(a_{j}-1)}{Z_{F_n}} W_{F_n}^{(j)}.
\end{align*}
the coefficient of $W_N^{(k)}$ ($1 \leq k \leq r-1$) is
\begin{align*}
w_{k+1}a_k + (a_k - 1) \sum_{i=1}^k w_i &= w_{k+1}a_k + (a_k - 1) a_k a_{k+1} \cdots a_{r-1} \\
&= a_k \{ w_{k+1} - (1-a_k) a_{k+1} a_{k+2} \cdots a_{r-1} \} = 0.
\end{align*}
Thus, we have
\begin{align*}
    \sum_{i=2}^r \frac{w_i a_{i-1}}{Z_{F_n}} W_{F_n}^{(i-1)} + \sum_{i=1}^{r-1} \sum_{j=i}^{r-1} \frac{w_i(a_{j}-1)}{Z_{F_n}} W_{F_n}^{(j)}=0.
\end{align*}
Therefore,
\begin{align*}
    \sum_{i=1}^r w_i H_{\nu_{n}^{(i)}} (\mathcal{U}^{(i)}_{F_n}) + w_1 \int_{X_1} S_{F_n} f d\nu_n=\log Z_{F_n}.
\end{align*}
Let $\mu\in \mathcal{M}^G(X)$ be the $G$-invariant probability measure on $X_1$ satisfying $(1),(2)$ above. Recall that
\begin{align*}
    Q^{\boldsymbol{a}}(G,f,\mathcal{U}^{(1)},\cdots,\mathcal{U}^{(r)})&=\lim_{n\to \infty}\displaystyle\frac{\log Q^{\boldsymbol{a}}(X_r, \spa f,\spa F_n, \spa \mathcal{U}^{(r)})}{|F_n|}\\&= \lim_{n\to \infty}\frac{1}{|F_n|}\log Z_{F_n}\\&= \lim_{n\to \infty}\frac{1}{|F_n|} \left(\sum_{i=1}^r w_i H_{\nu_{n}^{(i)}} (\mathcal{U}^{(i)}_{F_n}) + w_1 \int_{X_1} S_{F_n} f d\nu_n\right).
\end{align*}

    By Lemma \ref{lem2.13}, for any $E\in \mathcal{F}(G)$ we have
 \begin{align*}
    \frac{1}{|F_n|}H_{\nu_{n}^{(i)}} (\mathcal{U}^{(i)}_{F_n})&\leq \sum_{g\in F_n}\frac{1}{|F_n|}\left(\frac{1}{|E|}H_{\nu_{n}^{(i)}}(\mathcal{U}^{(i)}_{Eg})+|F\backslash\{g\in G:E^{-1}g\subset F_n\}|\cdot\log|\mathcal{U}^{(i)}|\right)\\&\leq\frac{1}{|E|}\left(\frac{1}{|F_n|}\sum_{g\in F_n}H_{\nu_{n}^{(i)}}(\mathcal{U}^{(i)}_{Eg})\right)+\frac{B(F_n,E^{-1}\cup\{e_G\})}{|F_n|}\cdot\log|\mathcal{U}^{(i)}|
 \end{align*}
 \begin{remark}
  Since $\varphi(x)=-x\log x$ is concave on $[0,1]$. For any measurable set $A\subset X$, we have

\begin{equation*}
    \varphi\Big(\displaystyle\frac{1}{|F_n|}\sum_{g\in F_n}g\nu^{(i)}_n(A)\Big)\geq\displaystyle\frac{1}{|F_n|}\sum_{g\in F_n}\varphi\Big(g\nu^{(i)}_n(A)\Big).
\end{equation*}
\end{remark}
Thus,
\[\frac{1}{|E|}\left(\frac{1}{|F_n|}\sum_{g\in F_n}H_{\nu_{n}^{(i)}}(\mathcal{U}^{(i)}_{Eg})\right)\leq \frac{1}{|E|}H_{\mu_n^{(i)}}(\mathcal{U}^{(i)}_E).\]
Therefore,
\begin{align*}
    &\lim_{n\to \infty}\frac{1}{|F_n|} \left(\sum_{i=1}^r w_i H_{\nu_{n}^{(i)}} (\mathcal{U}^{(i)}_{F_n}) + w_1 \int_{X_1} S_{F_n} f d\nu_n\right)\\&\leq\lim_{n\to \infty}\sum_{i=1}^{r}\frac{w_i}{|E|}\left(H_{\mu_n^{(i)}}(\mathcal{U}^{(i)}_E)+\frac{B(F_n,E^{-1}\cup\{e_G\})}{|F_n|}\cdot\log|\mathcal{U}^{(i)}|\right)+w_1\int_{X_1}fd\mu.
\end{align*}
Let $n\to \infty$, since $\mathcal{U}^{(i)}$ are clopen, we get
\begin{align*}
    Q^{\boldsymbol{a}}(G,f,\mathcal{U}^{(1)},\cdots,\mathcal{U}^{(r)})\leq \sum_{i=1}^rw_ih_{\mu_i}(\mathcal{U}^{(i)},G)+\omega_1\int_X fd\mu.
\end{align*}
\end{proof}

Therefore, by Remark \ref{rmk2.9} (2) for any clopen covers $\mathcal{U}^{(i)}\in \mathcal{C}_{X_i}^o$ $(i=1,\cdots,r)$, there is a $\mu\in \mathcal{M}^G(X_1)$ such that
    \begin{align*}
       Q^{\boldsymbol{a}}(G,f,\mathcal{U}^{(1)},\cdots,\mathcal{U}^{(r)})\leq\sum_{i=1}^rw_ih_{\mu_i}(X_i,G)+\omega_1\int_X fd\mu
    \end{align*}
    for zero-dimensional systems $(X_i,G)$ $(i=1,\cdots,r)$. Than by Claim \ref{cl4.3}, there is a $\mu\in \mathcal{M}^G(X_1)$ such that
    \[P^{\boldsymbol{a}}(G,f,\boldsymbol{\pi})\leq \sum_{i=1}^rw_ih_{\mu_i}(X_i,G)+\omega_1\int_X fd\mu.\]
\end{proof}
\begin{remark}
Let $(X_i, \spa G)$ ($i=1, \spa 2, \spa \ldots, \spa r$), $f\in C(X_1,\mathbb{R})$ and factor maps $\pi_i: X_i \rightarrow X_{i+1} \hspace{5pt} (i=1, \spa 2, \spa \cdots \spa , \spa r-1)$. From Proposition \ref{proposition: zero-dimensional trick}, there are zero-dimensional $G$-systems $(Z_i, \spa G)$ ($i=1, \spa 2, \spa \ldots, \spa r$), $g\in C(Z_1,\mathbb{R})$ and factor maps $\Pi_i: Z_i \rightarrow Z_{i+1} \hspace{5pt} (i=1, \spa 2, \spa \cdots \spa , \spa r-1)$ such that
\[ P^{\boldsymbol{a}}_{\mathrm{var}}(f, G, \boldsymbol{\pi}) \geq P^{\boldsymbol{a}}_{\mathrm{var}}(g, G, \boldsymbol{\Pi}) \]
and
\[ P^{\boldsymbol{a}}(f, G, \boldsymbol{\pi}) \leq P^{\boldsymbol{a}}(g, G, \boldsymbol{\Pi}). \]
Then by Theorem \ref{thm4.1}, we have
\[ P^{\boldsymbol{a}}(f, G, \boldsymbol{\pi}) \leq P^{\boldsymbol{a}}(g, G, \boldsymbol{\Pi})\leq P^{\boldsymbol{a}}_{\mathrm{var}}(g, G, \boldsymbol{\Pi})\leq P^{\boldsymbol{a}}_{\mathrm{var}}(f, G, \boldsymbol{\pi}). \]
Therefore, Theorem \ref{thm3.1} is obtained.
\end{remark}

\section{Weighted pressure determines weighted measure-theoretic entropy}
In this section, we will give a simple application of the weighted variational principle for pressure. We investigate how the weighted pressure determines the weighted measure-theoretic entropy. At the beginning, we need the following lemma, one can refer to \cite{ds} for a proof.
\begin{lemma}\label{l1}
Let $V$ be a locally convex linear topological space and $K_1, K_2$ be closed and convex with $K_1\cap K_2=\emptyset$. If additionally $K_1$ is compact, then there is a continuous real-valued linear functional $F$ on $V$ such that $F(X)<F(y)$ for all $x\in K_1$ and $y\in K_2$.
\end{lemma}
\begin{theorem}
Let $\mu_0\in \mathcal{M}^G(X_1)$. Suppose that $h_{top}^{\boldsymbol{a}}(X_1,G)<+\infty$ and the entropy maps $\theta\in \mathcal{M}^G(X_i)\mapsto h_{\theta}(X_i,G)$, $i=1,\cdots,r$ are upper semi-continuous at $\mu_0$. Then
\[
h_{\mu_0}^{\boldsymbol{a}}(X_1,G)=\inf\left\{P^{\boldsymbol{a}}(G,f)-\int_{X_1}fd\mu_0|f\in C(X_1)\right\}.
\]
\end{theorem}
\begin{proof}
By the variational principle in Theorem \ref{main}, we have
\[
h_{\mu_0}^{\boldsymbol{a}}(X_1,G)\leq\inf\left\{P^{\boldsymbol{a}}(G,f)-\int_{X_1}fd\mu_0|f\in C(X_1)\right\}.
\]
We now prove the opposite inequality. Let $b>h_{\mu_0}^{\boldsymbol{a}}(X_1,G)$ and let
\[
C=\{(\mu,t)\in \mathcal{M}^G(X_1)\times \mathbb{R}|0\leq t\leq h_{\mu}^{\boldsymbol{a}}(X_1,G)\}.
\]
\textbf{Claim:} $C$ is a convex set.

In fact, let $(\mu_1,t_1),(\mu_2,t_2)\in C$ and $p\in[0,1]$. Then we have $0\leq t_1\leq h_{\mu_1}^{\boldsymbol{a}}(X_1,G)$ and $0\leq t_2\leq h_{\mu_2}^{\boldsymbol{a}}(X_1,G)$. Since the entropy map $\mu\in \mathcal{M}^G(X_1)\mapsto h^{\boldsymbol{a}}_{\mu}(X_1,G)$ is affine, we have
\[
h_{p\mu_1+(1-p)\mu_2}^{\boldsymbol{a}}(X_1,G)=ph_{\mu_1}^{\boldsymbol{a}}(X_1,G)+(1-p)h_{\mu_2}^{\boldsymbol{a}}(X_1,G)\geq pt_1+(1-p)t_2\geq0.
\]
Then $p(\mu_1,t_1)+(1-p)(\mu_2,t_3)\in C$ and the claim is proved.

Let $C(X_1)^*$ be the dual space of $C(X_1)$ endowed with the weak$^*$ topology and consider $C$ as a subset of $C(X_1)^*\times\mathbb{R}$. Then by the upper semi-continuity of the entropy map at $\mu_0$, we have $(\mu_0,b)\notin \overline{C}$. Applying Lemma \ref{l1} to the disjoint convex sets $\overline{C}$ and $(\mu_0,b)$, we can get a continuous real-valued linear functional $F$ on $C(X_1)^*\times\mathbb{R}$ such that $F((\mu,t))<F((\mu_0,b))$ for all $(\mu,t)\in \overline{C}$. By the weak$^*$ topology on $C(X_1)^*$, we have that there exist $f\in C(X_1)$ and $d\in\mathbb{R}$ such that
\[
F((\mu,t))=\int_{X_1}fd_{\mu}+dt.
\]
It follows that $\int_{X_1}fd_{\mu}+dt<\int_{X_1}fd_{\mu_0}+db$. In particular, by the construction of $C$, we have
\begin{equation}\label{eq5.1}
\int_{X_1}fd_{\mu}+dh_{\mu}^{\boldsymbol{a}}(X_1,G)<\int_{X_1}fd_{\mu_0}+db
\end{equation}
for all $(\mu,t)\in \overline{C}$.
Let $\mu=\mu_0$ in \eqref{eq5.1}. Then by the assumption that $b>h_{\mu_0}^{\boldsymbol{a}}(X_1,G)$, we have $d>0$. Hence,
\begin{equation*}
\int_{X_1}\frac{f}{d}d_{\mu}+h_{\mu}^{\boldsymbol{a}}(X_1,G)<\int_{X_1}\frac{f}{d}d_{\mu_0}+b,\text{ for all $\mu\in \mathcal{M}^G(X_1)$}.
\end{equation*}
Therefore, by the variational principle in Theorem \ref{main}, we have
\[
P^{\boldsymbol{a}}(G,\frac{f}{d})\leq \int_{X_1}\frac{f}{d}d_{\mu_0}+b.
\]
Rearranging gives
\[
b\geq P^{\boldsymbol{a}}(G,\frac{f}{d})-\int_{X_1}\frac{f}{d}d_{\mu_0}\geq \inf\left\{P^{\boldsymbol{a}}(G,g)-\int_{X_1}gd\mu_0|g\in C(X_1)\right\}.
\]
Then $h_{\mu_0}^{\boldsymbol{a}}(X_1,G)\geq \inf\left\{P^{\boldsymbol{a}}(G,g)-\int_{X_1}gd\mu_0|g\in C(X_1)\right\}.$
\end{proof}
\section*{Acknowledgement}
The authors are partly supported by NNSF of China (Grant No. 12201120).

\vspace{0.5cm}

\address{Department of Mathematics, Nanjing University, Nanjing 210093, People's Republic of China}

\textit{E-mail}: \texttt{yzy199707@gmail.com}

\address{School of Mathematics and Statistics, Fuzhou University, Fuzhou 350116, People's Republic of China}

\textit{E-mail}: \texttt{xzb2020@fzu.edu.cn}

\end{document}